\documentclass{amsart}

\newtheorem{theorem}[equation]{Theorem}
\newtheorem{prop}[equation]{Proposition}
\newtheorem{lemma}[equation]{Lemma}

\theoremstyle{remark}
\newtheorem{remark}[equation]{Remark}
\theoremstyle{definition}

\numberwithin{equation}{subsection}

\renewcommand{\qed}{\hspace*{\fill} \setlength{\unitlength}{1mm}
\begin{picture}(2.5,2.5)
      \put(0,0){\framebox(2.5,2.5){}}
  \end{picture}
\setlength{\unitlength}{1pt}}

\newcommand{\diam}{\operatorname{diam}}

\newcommand{\inj}{\operatorname{inj}}
\newcommand{\Lsp}{\operatorname{Lsp}}

\newcommand{\Tr}{\operatorname{Tr}}

\newcommand{\reals}{{\mathbb R}}

\newcommand{\cP}{{\mathcal{P}}}
\newcommand{\cR}{{\mathcal{R}}}

\newcommand{\cA}{{\mathcal{A}}}

\newcommand{\cH}{{\mathcal{H}}}
\begin{document}
\title[Weyl's law on negatively curved surfaces]
{A lower bound for the remainder in Weyl's law on negatively curved
surfaces} 

\author[D. Jakobson]{Dmitry Jakobson}
\address{Department of Mathematics and
Statistics, McGill University, 805 Sherbrooke Str. West,
Montr\'eal QC H3A 2K6, Ca\-na\-da.}
\email{jakobson@math.mcgill.ca}

\author[I. Polterovich]{Iosif Polterovich}
\address{D\'e\-par\-te\-ment de math\'ematiques et de
sta\-tistique, Univer\-sit\'e de Mont\-r\'eal CP 6128 succ
Centre-Ville, Mont\-r\'eal QC  H3C 3J7, Canada.}
\email{iossif@dms.umontreal.ca}

\author[J. Toth]{John A. Toth}
\address{Department of Mathematics and
Statistics, McGill University, 805 Sherbrooke Str. West,
Montr\'eal QC H3A 2K6, Ca\-na\-da.}
\email{jtoth@math.mcgill.ca}



\thanks{D.J. and J.T. are supported by NSERC, FQRNT and Dawson Fellowships.
I.P. is supported by NSERC and FQRNT}

\begin{abstract}
We obtain an estimate from below for the remainder in Weyl's law on
negatively curved surfaces. In the constant curvature case, such a
bound was proved  independently by Hejhal and Randol in 1976 using
the Selberg zeta function techniques. Our approach works in
arbitrary negative curvature, and is based on wave trace asymptotics
for long times, equidistribution of closed geodesics and small-scale
microlocalization.
\end{abstract}
\maketitle


\section{Introduction and main results}\label{sec:lowglobal}
\subsection{Weyl's law}
Let $X$ be a compact negatively curved surface of area $A$ with  the
Riemannian metric $\{g_{ij}\}$. We assume that the Gaussian
curvature $K(x)$  satisfies
\begin{equation}
\label{curvature}
-K_1^2 \le K(x) \le -K_2^2
\end{equation}
at every point $x \in X$. Let $\Delta$ be the Laplacian on $X$ with
the eigenvalues $0=\lambda_0<\lambda_1 \leq \lambda_2 \leq \dots$
and the corresponding orthonormal basis $\{\phi_i\}$ of
eigenfunctions: $\Delta \phi_i=\lambda_i\phi_i$. Let
$N(\lambda)=\#\{\sqrt{\lambda_i} \le \lambda\}$ be the eigenvalue
counting function. The asymptotic behavior of $N(\lambda)$ is given
by {\it Weyl's law} (\cite{H}):
\begin{equation}\label{Weyl}\begin{aligned}
N(\lambda)&=\frac{A}{4\pi}\lambda^2 + R(\lambda),
\qquad R(\lambda)=O(\lambda).
\end{aligned}
\end{equation}
It shown in \cite{Berard} that on a nonpositively curved surface
$R(\lambda)=O\left(\frac{\lambda}{\ln \lambda}\right)$. In the
present paper we study lower bounds for $R(\lambda)$. As in
\cite{JP}, one of our tools is thermodynamic formalism for
hyperbolic flows.


\subsection{Thermodynamic formalism}
Let $G^t$ be the geodesic flow on the unit tangent bundle $SX$ and
let $E_\xi^u$ be the (one-dimensional) unstable subspace for $G^t$,
$\xi \in SX$. The {\it Sinai-Ruelle-Bowen potential} is a H\"older
continuous function $\cH:SX \to\reals$ which   for any $\xi\in SX$
is defined by the formula (see \cite{BR}, \cite{Sinai2})
\begin{equation}\label{SRB}
\cH(\xi)=\left.\frac{d}{dt}\right|_{t=0}\ln\det dG^t|_{E_\xi^u},
\end{equation}
For any continuous function $f:SX \to \reals$ one can define the
{\it topological pressure}
\begin{equation}\label{pressure}
P(f)=\sup_\mu \left(h_\mu + \int f d\mu \right),
\end{equation}
where the supremum is taken over all $G^t$-invariant measures $\mu$
and $h_\mu$ denotes the {\it measure-theoretical entropy} of the
geodesic flow (see \cite{Bowen}).
 In particular $P(0)=h$,
where $h$ is the {\it topological entropy} of the flow. It is
well-known that for the Sinai-Ruelle-Bowen potential $P(-\cH)=0$ and
the corresponding equilibrium measure (i.e. the measure on which the
supremum is attained) is the Liouville measure $\mu_L$ on the unit
tangent bundle:
\begin{equation}\label{ent}
h_{\mu_L} = \int_{SX} \cH d\mu_L
\end{equation}

\subsection{Main result}
Recall that $f_1(\lambda)=\Omega(f_2(\lambda))$  for a function
$f_1$ and a positive function $f_2$ means that $\limsup_{\lambda \to
\infty} |f_1(\lambda)|/f_2(\lambda)>0$.
\begin{theorem}\label{main}
Let $X$ be a compact surface of negative curvature. Then
\begin{equation}
\label{bound} R(\lambda)=\Omega
\left((\ln\lambda)^{\frac{P(-\cH/2)}{h}-\varepsilon}\right) \,\,\,
\forall \, \varepsilon >0.
\end{equation}
\end{theorem}
As was shown in \cite[section 1.4]{JP}, the power of the logarithm
in \eqref{bound} is always positive:
$$
\frac{P\left(-\cH/2\right)}{h} \ge \frac{K_2}{2K_1}>0.
$$
Moreover, if the curvature is constant, the bound \eqref{bound}
reads
\begin{equation} \label{constcurv}
R(\lambda)=\Omega \left((\ln
\lambda)^{\frac{1}{2}-\varepsilon}\right)
\end{equation}
for any $\varepsilon >0$. The estimate \eqref{constcurv} was
obtained in  \cite{Randol} and, in a slightly stronger form, in
\cite[section 17]{Hej}, using the Selberg zeta function techniques.
Our approach, based on the Duistermaat-Guillemin wave trace formula,
thermodynamic formalism and semiclassical analysis, allows us to
treat the {\it variable} curvature case as well.

Theorem \ref{main} agrees with a ``folklore'' conjecture that on a
generic negatively curved surface
\begin{equation}
\label{folklore} R(\lambda)=O(\lambda^\varepsilon)  \,\,\,\,\,
\forall \, \varepsilon >0.
\end{equation}
Genericity is important, since on arithmetic surfaces corresponding
to quaternionic lattices one can prove a much better lower bound
$R(\lambda)= \Omega \left( \frac{\sqrt{\lambda}}{\ln
\lambda}\right)$ (see \cite{Hej}).

\begin{remark}One may compare Theorem \ref{main} with
the lower bound for the {\it pointwise} error term obtained in
\cite{JP}:
\begin{equation}
\label{pointwise} R_x(\lambda)=\sum_{\sqrt{\lambda_i}\le \lambda}
|\phi_i(x)|^2 - \frac{\lambda^2}{4\pi} =
\Omega\left(\sqrt{\lambda}\,\,
(\ln\lambda)^{\frac{P(-\cH/2)}{h}-\varepsilon}\right) \,\,\,\,\,
\forall \, \varepsilon >0. \end{equation}
Estimates \eqref{bound} and \eqref{pointwise} are independent, and,
in particular, \eqref{bound}
can not be deduced from \eqref{pointwise}. 
Indeed, cancellations may occur when $R_x(\lambda)$ is integrated
over a negatively curved surface (for instance, according to the
conjecture \eqref{folklore}, $\sqrt{\lambda}$ should cancel out in
the generic case). 
Also, the sequence of $\lambda$--s yielding
the $\Omega$--bound \eqref{pointwise} depends on the lengths of the
geodesic loops at
$x$ (see \cite[section 5.2]{JP}), and hence for each point $x$ such a sequence is apriori different. 
\end{remark}


\subsection{Wave trace asymptotics for long times}
\label{sec:DG} Consider the spectral distribution
\begin{equation}\label{waveker1}
\Tr \, e(t)\ =\ \sum_{i=0}^\infty \cos(\sqrt{\lambda_i}t)
\end{equation}
which is the even part of the wave trace on $X$.

To prove Theorem \ref{main} we use a modification of the
Duisterma\-at-Guillemin asymptotic formula for the wave trace
\cite{DG}. Originally, this formula captures the contribution of a
single closed geodesic, however  for the proof of Theorem \ref{main}
we need to take into account the contributions of {\it all} closed
geodesics of length $T_0 < L \le T(\lambda)$, where $T_{0} >0$ is
some constant (see Lemma \ref{dynlemma}), and $T(\lambda) \to
\infty$ at an appropriate rate as the spectral parameter $\lambda
\to \infty$.

Let $\chi(t,T)$ be a cut-off function
\begin{equation} \label{cutoff0}
\chi(t,T) =  ( 1 - \psi(t) ) \, \hat{\rho}\left(\frac{t}{T}\right),
\end{equation}
where $\rho \in {\mathcal S}(\reals)$  is  an even, non-negative
Schwartz function such that ${\rm supp}\, \hat{\rho} \subset
[-1,+1]$, and $\psi(t) \in C^{\infty}_{0}(\mathbb{R})$ with $\psi(t)
\equiv 1$ when $ t \in [-T_{0}, T_{0}]$ and $\psi(t) \equiv 0$ when
$|t| \geq 2T_{0}$.

The long-time version of the Duistermaat-Guillemin trace formula is
given by
\begin{theorem}
\label{tk2} Let $T(\lambda) \rightarrow \infty$ as $\lambda
\rightarrow \infty$  with $T(\lambda) \leq\epsilon \ln \lambda$ for
$\epsilon >0$  small enough, and let $L_\gamma$ (resp.
$L_\gamma^\sharp$) denote the length (resp. the primitive period) of
a periodic geodesic $\gamma$. Then the asymptotics of the smoothed
Fourier transform of the wave trace is given by
\begin{equation}\label{tqq}
\int_{-\infty}^{\infty} \Tr \, e(t) \chi(t,T) \cos \lambda t \, dt
=\sum_{L_{\gamma} \in \Lsp, L_{\gamma}\le T(\lambda)}
\frac{L_{\gamma}^\sharp \cos (\lambda L_{\gamma})\,\chi(L_{\gamma},
T)} {\sqrt{|\det(I-\cP_{\gamma})|}} + {\mathcal O}\left(1 \right).
\end{equation}
Here, $\cP_{\gamma}$ is the linearized  Poincar\'{e} map corresponding
to $\gamma$.
\end{theorem}

The proof of Theorem \ref{tk2}  is the most technically difficult
part of the paper. We use semiclassical microlocal analysis on small
scales, see section 3. In order to avoid the ``accumulation'' of
singularities in the wave trace and to make the stationary phase
method work, we separate the contributions of each closed geodesic
using Lemma \ref{dynlemma}. Although the idea of working to
(suitably scaled) Ehrenfest times $T(\lambda) \sim \epsilon \ln
\lambda$ is well-established (see \cite{Berard, Zel2, Faure}),
rigorously separating out large exponential sums with
$T(\lambda)$-terms from the wave-trace using small-scale
$\hbar$-microlocalization appears to be a novel approach to
estimating remainders in the negatively-curved case.

The main term in the asymptotics proved in Theorem \ref{tk2} is
given by the sum of the principal wave invariants at each closed
geodesic. As follows from Lemma \ref{above:weak1}, to prove Theorem
\ref{main} it is sufficient to show that this sum grows at the rate
given by  \eqref{bound}. First, we  use the Parry-Pollicott
equidistribution result \cite{PP} to calculate the asymptotics of
the sum not taking into account its oscillatory nature. Then we deal
with the oscillations in the wave invariants (the difficulty is that
oscillating terms may cancel out in the sum) using a ``straightening
the phases'' argument based on the Dirichlet box principle (cf.
\cite{JP}).

\subsection{Lower bound for $R(\lambda)$ in higher dimensions}
Lower bounds for the error term in Weyl's law on
higher-dimensional manifolds are in fact much simpler than on
surfaces. The reason is that the contribution from the
``singularity at zero'' to the remainder dominates the
contribution from the periodic geodesics.

We write $f(\lambda) \gg \lambda^k$ for a function $f(\lambda)$ if
there exist a constant $c_0>0$ and a number $\lambda_0$, such that
$f(\lambda)>c_0 \lambda^k$ for any $\lambda
> \lambda_0$. The following $L^1$-estimate  (which is
stronger than an $\Omega$-bound) holds for the remainder in Weyl's
law on a rather general class of Riemannian manifolds containing manifolds
of negative curvature.
\begin{theorem}\label{Riesz}
Let $X$ be a manifold of dimension $n\ge 3$, such that $\int_X
\tau  \neq 0$, where $\tau$ is the scalar curvature. Then,
\begin{equation}\label{bound3}
\frac{1}{\lambda}\int_0^{\lambda} |R(t)|dt \gg \lambda^{n-2},
\end{equation}
\end{theorem}
Theorem \ref{Riesz}  is proved using the asymptotics of the Riesz
means (see \cite{Saf2}) in section \ref{sec:Riesz}. One can also
prove Theorem \ref{Riesz} using the standard  $t \rightarrow 0^{+}$
heat trace asymptotics, see \cite[\S 2.1]{JP}. In order to get more
refined information about $R(\lambda)$ on negatively curved
manifolds of higher dimension it is natural to study the {\it
oscillatory} error term, $R^{osc}(\lambda)$, see section
\ref{oscerr}.


\section{Two auxiliary lemmas}
\subsection{Smoothed Fourier transform of the wave trace}
In the notations of section \ref{sec:DG} let
\begin{equation}\label{klamt}
k(\lambda,T)=\int_{-\infty}^\infty\frac{\hat{\rho}(t/T)}{T}
\cos(\lambda t) \Tr \, e(t) dt,
\end{equation}
where, $T=T(\lambda)$ will be chosen appropriately later on.
Substituting  \eqref{waveker1} into \eqref{klamt} we obtain
\begin{equation}\label{KH}
k(\lambda, T)=\sum_{i=0}^\infty
H_{\lambda,T}(\sqrt{\lambda_i}),
\end{equation}
where, for $r \geq 0,$
\begin{multline}\label{defH}
H_{\lambda,T}(r)=\int_{-\infty}^\infty\frac{\hat{\rho}(t/T)}{T}
\cos(\lambda t) \cos(r t) dt=\frac{ \rho(T(\lambda -r )) +
\rho(T(\lambda+r))}{2}= \\ \frac{ \rho(T(\lambda -r ) )}{2} +
{\mathcal O}(\lambda^{-\infty}),
\end{multline}
Here,
$$
\widehat{\rho}(s)= \frac{1}{2\pi}\int_{-\infty}^{\infty}
e^{-is\zeta}\rho(\zeta)d\zeta
$$
is the Fourier transform of $\rho$ and the  ${\mathcal
O}(\lambda^{-\infty})$-error in (\ref{defH}) follows from the fact
that  $\rho (\lambda+ r) = {\mathcal O}(\lambda^{-\infty})$
uniformly for $r \geq 0$, since $\rho \in {\mathcal S}({\mathbb
R})$. Replacing the sum in \eqref{KH} by an integral, we get the
following representation of $k(\lambda,T)$:
\begin{equation}\label{klamt2}
k(\lambda,T) = \int_0^\infty
\frac{\rho(T(\lambda-r))+\rho(T(\lambda+r))}{2}dN(r)= \int_0^\infty
H_{\lambda,T}(r) dN(r)
\end{equation}
Formula \eqref{klamt2} plays a key role in our analysis.
We shall also use the following notation:
\begin{equation}\label{klamt3}
\int_0^\infty H_{\lambda,T}(r)
dR(r)= \kappa(\lambda,T)
\end{equation}
Note that the contribution of the main term in Weyl's law has been
subtracted from $k(\lambda,T)$ to get $\kappa(\lambda,T)$.


We shall use the following:
\begin{lemma}
\label{above:weak1}  Let $R(\lambda)=o((\ln \lambda)^b)$, $b>0$.
Then $\kappa(\lambda, T)=o((\ln \lambda)^b)$ uniformly in $T$ for
$T>T_0$, where $T_0$ is an arbitrary positive number.
\end{lemma}
\noindent{\bf Proof.}  By the assumption of the lemma, for any
$\epsilon>0$, $R(\lambda)< \epsilon (1+\ln \lambda)^b$ for large
enough $\lambda$. Consider the left hand side of \eqref{klamt3}:
\begin{equation}
\label{LL} \int_0^\infty H_{\lambda,T}(r) dR (r).
\end{equation}
Taking into account \eqref{defH} and  integrating \eqref{LL} by
parts we obtain
\begin{equation}
\label{byparts} \kappa(\lambda, T) \le \frac{\epsilon
T}{2}\int_0^\infty |\rho'(T(r-\lambda))| (\ln (1+r))^b dr
+\frac{\epsilon T}{2}\int_0^\infty |\rho'(T(r+\lambda))| (\ln
(1+r))^b dr.
\end{equation}
Since $\rho'$ is Schwartz class, the second term of \eqref{byparts}
is $O(1)$.  Changing variables in the first term of \eqref{byparts},
we obtain
\begin{multline*}
\frac{\epsilon T}{2}\int_{-\lambda}^{\infty} |\rho'(Ts)|
(\ln(1+\lambda+s))^b ds=\\\frac{\epsilon (\ln\lambda)^b}{2}
\int_{-\lambda T}^{\infty} |\rho'(u)|
\left(1+\frac{\ln(1+\frac{u+T}{\lambda T})}{\ln\lambda}\right)^bdu
\le C \epsilon (\ln\lambda)^b
\end{multline*}
for some constant $C>0$, where the last inequality again follows
from the fact that $\rho'$ is Schwartz class. Clearly, the constant
$C$ can be chosen uniformly in $T$ for $T>T_0>0$. Since $\epsilon$
can be taken arbitrarily small, we get $\kappa(\lambda,
T)=o((\ln\lambda)^b)$, and this completes the proof of the lemma.
\qed
\begin{remark} Lemma \ref{above:weak1} is proved similarly to the results
of \cite[section 2.2]{JP}, see also \cite{K}. One may also compare
it to \cite[Proposition 3.1]{Sar}. Sarnak's argument gives the lower
bound $R(\lambda)=\Omega(\sqrt{\lambda})$ for the Weyl error on a
surface under the assumption  that  the geodesic flow $G^t$ has a
fixed point set of dimension two for some $t>0$. This condition
holds, for example, when $G^t$ is completely integrable, but it is
not satisfied on a negatively curved surface.
\end{remark}

\subsection{Separation of periodic orbits}\label{sec:dynlemma}
In this section we shall prove the following dynamical lemma:
\begin{lemma}\label{dynlemma}
Let $X$ be a negatively curved surface and let $
\varOmega(\gamma,\varepsilon)$ denote the $\varepsilon$-neighborhood
of a geodesic $\gamma$ in $SX$ with respect to the Sasaki metric.
Then there exist positive constants $T_0, B$ and  $\delta'$
(depending only on the injectivity radius $\inj(X)$ and the lower
curvature bound $K_{1}$) such that for any $T>T_0$ the sets
$\varOmega(\gamma,e^{-BT})$ are disjoint for all pairs of closed
geodesics $\gamma$ on $X$ with length $L_\gamma\in [T-\delta',T].$
\end{lemma}

Note that since there are exponentially many closed geodesics on $X$
of length $L_\gamma\in [T-\delta',T]$, disjoint neighborhoods have
to be of exponentially small size.

\begin{remark}
Here and further on we write $\varOmega(Y,d)$ for the
$d$-neighborhood of the set $Y$. It should not be confused with the
$\Omega$ notation for the lower bounds.
\end{remark}

\smallskip

\noindent {\bf Proof.} Choose $B\geq 2 K_1$, where $K_1^2$ is the
curvature bound and $K_1$ is an upper bound for Lyapunov exponents,
cf. \eqref{eigunstable}.  Also, choose $\delta'<{\rm inj}(X)/3$, and
let $T_0$ be such that $2e^{-K_1 T_0}<\delta'$. Assume for
contradiction that there exist two closed geodesics $\gamma_1$ and
$\gamma_2$ with $T-\delta' \leq L_{\gamma_1} \leq L_{\gamma_2}\leq
T$ such that the corresponding neighborhoods intersect, and that the
geodesics are not inverses of each other (note that due to the
restriction $L_\gamma\in [T-\delta',T]$ in the conditions of the
lemma, the geodesics cannot be integer multiples of each other
unless they are inverses).

Denote the geodesics on $X$ by $\gamma_j(t),0\leq t\leq
L_{\gamma_j}$, and their lifts to $SX$ by $\tilde
\gamma_j(t)=(\gamma_j(t),\gamma_j'(t))$;  $\gamma_j(t)\subseteq X$
is sometimes called a {\em footprint} of $\tilde \gamma_j(t)$.
Without loss of generality we may assume that
$$
{\rm dist}_{SX}(\tilde \gamma_2(0),\tilde \gamma_1(0))\leq
2e^{-2K_1T}.
$$
Since for any $0\le t \le L_{\gamma_2}$
$$
{\rm dist}_{SX}(\tilde \gamma_2(t),\tilde \gamma_1(t))= {\rm
dist}_{SX}(G^t\tilde \gamma_2(0),G^t \tilde \gamma_1(0)) \le
2e^{-2K_1T} e^{K_1t} \le 2e^{-K_1T},
$$
and hence
\begin{equation}\label{dist1:base}
{\rm dist}_X(\gamma_2(t), \gamma_1(t))\le 2e^{-K_1T}.
\end{equation}
In other words, the entire geodesic $\gamma_2$ lies in the
$2e^{-K_1T}$-neighborhood of the geodesic $\gamma_1$ and vice versa.

For convenience, we reparametrise $\gamma_1$ and define
$$
\beta_1(s):=\gamma_1(L_1s/L_2),\qquad 0\leq s\leq L_2.
$$
By triangle inequality and the definition of $\beta_1$,
\begin{equation}\label{dist2:base}
\begin{aligned}
\;& {\rm dist}(\gamma_2(t), \beta_1(t))\leq {\rm dist}(\gamma_2(t),
\gamma_1(t)) + {\rm dist}(\gamma_1(t), \beta_1(t)) \leq\\
\;& 2e^{-K_1t} + t\left(1- \frac{L_1}{L_2} \right)\leq
2e^{-K_1T}+\delta'<\frac{2\cdot{\rm inj}(X)}{3}.
\end{aligned}
\end{equation}

Accordingly, for any $0\leq t\leq L_2$ there exists a unique
shortest geodesic $\alpha_t(s)$ in $X$ connecting $\gamma_2(t)$ and
$\beta_1(t)$. We shall choose the parameter $s\in[0,1]$ so that
$\alpha_t(0)=\gamma_2(t)$ and $\alpha_t(1)=\beta_1(t)$.

Define the mapping $\Phi(t,s):[0,L_2]\times[0,1]\to X$ by the
formula
$$
\Phi(t,s)=\alpha_t(s).
$$
We claim that $\Phi$ defines a homotopy between $\gamma_2(t)$ and
$\beta_1(t)$.  Indeed, $\Phi(t,0)=\gamma_2(t)$,
$\Phi(t,1)=\beta_1(t)$.  Moreover, since both $\gamma_2$ and
$\beta_1$ have period $L_2$, we have
$$
\alpha_0(s)=\alpha_{L_2}(s), \ \forall\; s\in[0,1],
$$
and so $\Phi(\cdot,s)$ is a closed curve in $X$.  Finally,
$\Phi(t,s)$ is continuous since the function ${\rm
dist}(\gamma_2(t),\beta_1(t))$ is a continuous function of $t$.

On the other hand, $\beta_1$ is just a reparametrization of
$\gamma_1$, hence $\gamma_1$ and $\gamma_2$ lie in the same free
homotopy class, contradicting the fact that on a negatively-curved
surface there is a unique closed geodesic in each free homotopy
class. The contradiction completes the proof of the lemma. \qed

\begin{remark}\label{cot}
Lemma \ref{dynlemma} is proved for the tangent bundle. The tangent
and the cotangent bundles can be identified using the Riemannian
metric, and therefore Lemma \ref{dynlemma} holds in the cotangent
bundle as well. In this setting it will be used in the next section.
\end{remark}

\begin{remark}The proof of Lemma \ref{dynlemma}
generalizes verbatim to higher dimensions.
\end{remark}

\begin{remark}
In the literature on wave invariants \cite{DG, Don, Zelditch} it is
customary to choose a cut-off function in trace formulas in the
length spectrum so that only a single length of a closed geodesic is
contained in its support.  Since we go to
$T(\lambda)\sim\ln\ln\lambda$ times, it is not enough to localize
exclusively in the length spectrum.  Indeed, when we localize around
a geodesic(s) of length $L_i$, error terms in expressions like
\eqref{qq} are of the order $O(1/|L_{i+1}-L_i|)$, where
$\ldots<L_i<L_{i+1}<\ldots $ denote distinct lengths of closed
geodesics on $X$ (ignoring the multiplicity). Therefore, the error
would be large in the presence of ``near-multiplicities''. We can
not control the gaps in the length spectrum on a generic negatively
curved surface, and hence we have to localize in the phase space as
it is done in the next section.
\end{remark}

\section{Wave trace asymptotics and small-scale microlocalization}
\label{sec:mainproof} In this section we give a proof of  Theorem
\ref{tk2} which  is quite technical.  Let us note that assuming
Proposition \ref{k2}, the subsequent sections can be read
independently of section 3.

\subsection{Plan of the argument}

We choose a parameter $T =
T(\lambda) > 0$ satisfying:
\begin{itemize} \label{period bounds}
\item{ (i) \, $T(\lambda) \rightarrow \infty$ \, as \, $\lambda \rightarrow \infty$ }
\item{ (ii) \, $T(\lambda) \leq \epsilon \ln \lambda.$ }
\end{itemize}
Here, $\epsilon >0$ is a small constant that will be chosen later
on. For the  applications in this paper, it will only be necessary
to take $T(\lambda)$ of order $\ln \ln \lambda$, so we will not be
concerned with determining the best possible constant $\epsilon >0$
in (\ref{period bounds}).

To localize the contribution to the wave trace and to
$\kappa(\lambda,T)$ from a given closed geodesic $\gamma$ of length
$L_\gamma \leq T$ it is important to microlocalize the wave trace
$\Tr e(t)$ to a neighborhood of $\gamma$ and then sum over all the
different $\gamma$'s. The complication here is that we want to take
into account the contributions of {\em all} closed geodesics
$\gamma$ with $T_0<L_\gamma \leq T(\lambda).$  Since $T(\lambda)
\rightarrow \infty$ as $\lambda \rightarrow \infty$, the number of
these geodesics blows up. As a result, simply summing stationary
phase expansions for each of the $\gamma$'s (which is automatic when
$T = {\mathcal O}(1)$) is impossible when one needs to work with
such long period intervals. The way to deal with this is to
microlocalize on neighborhoods of the $\gamma$'s that shrink fast
enough as  $\lambda \rightarrow \infty$  (but not too fast) and then
split up the time interval $[T_{0}, T(\lambda)]$ into short
``windows'' of fixed size $\delta'>0$. In this context, it is
natural to work with semiclassical pseudodifferential and Fourier
integral operators. The crucial dynamical result we need here is
Lemma \ref{dynlemma}. Roughly speaking, this lemma says that there
exist neighborhoods $\varOmega(\gamma,e^{-BT})$ of $\gamma$ of size
$e^{-BT}$ in phase space $T^{*}X$ with the property that, for
appropriate geometric constant $\delta'>0$, no other periodic
geodesic with period in the window $[L_\gamma- \delta', L_\gamma +
\delta']$ intersects $\varOmega(\gamma,e^{-BT})$. Since $T(\lambda)
\leq \epsilon \ln \lambda,$ this clearly suggests microlocalizing
the trace to $e^{-B \epsilon \ln \lambda} =
\lambda^{-B\epsilon}$-neighborhoods of $\gamma$.

\subsection{A reformulation of Theorem \ref{tk2}}
For the purposes of the proof of Theorem \ref{main}, it is
convenient for us to reformulate Theorem \ref{tk2} in the following
equivalent form:
\begin{prop}
\label{k2} Let $T(\lambda) \rightarrow \infty$ as $\lambda
\rightarrow \infty$  with $T(\lambda) \leq\epsilon \ln \lambda$.
Then, for $\epsilon >0$  small enough,
\begin{equation}\label{qq}
\kappa(\lambda,T(\lambda))=\sum_{L_{\gamma} \in \Lsp, L_{\gamma}\le
T(\lambda)}  \frac{L_{\gamma}^\sharp \cos (\lambda
L_{\gamma})\,\chi\left(L_{\gamma}, T(\lambda) \right)}{ T(\lambda)
\sqrt{|\det(I-\cP_{\gamma})|}} + {\mathcal
O}\left(\frac{1}{T(\lambda)}\right).
\end{equation}
Here, $\cP_{\gamma}$ is the linearized  Poincar\'{e} map
corresponding to $\gamma$ and $\chi(t,T)$ is a cut-off function
defined by \eqref{cutoff0}.
\end{prop}
The proof of Proposition \ref{k2} is divided in several steps that
are carried out in  sections 3.4--3.8.

Equivalence of Theorem \ref{tk2} and Proposition \ref{k2}
immediately follows from the following
\begin{lemma}
\label{reform}
  Let $T(\lambda) \rightarrow \infty$ as $\lambda
\rightarrow \infty$.
Then
\begin{equation} \label{split1}
\kappa(\lambda,T) = \frac{1}{T(\lambda)}  \int_{-\infty}^{\infty}
\Tr \, e(t) \, \chi(t, T)  \cos ( \lambda t) \, dt + {\mathcal
O}\left(\frac{1}{T(\lambda)}\right).
\end{equation}
\end{lemma}
\noindent{\bf Proof.}  Combining \eqref{klamt}, \eqref{klamt2},
\eqref{klamt3} and using that $\rho$ is a Schwartz function we get
\begin{multline} \label{prelim1}
\kappa(\lambda,T) =  \int_{-\infty}^{\infty}
\frac{\hat{\rho}(t/T)}{T} \cos(\lambda t) \Tr \, e(t) dt -\\
\frac{A}{4\pi}\int_0^\infty r(\rho(T(\lambda-r))+\rho(T(\lambda+r)))
dr =\\ \int_{-\infty}^{\infty}\frac{\hat{\rho}(t/T)}{T} \cos(\lambda
t) \Tr \, e(t) dt - \frac{A}{4\pi}\int_0^\infty r
\rho(T(\lambda-r))dr + {\mathcal O}((\lambda T)^{-\infty})=
\\ \int_{-\infty}^{\infty} \frac{\hat{\rho}(t/T)}{T} \cos(\lambda t)
\Tr \, e(t) dt - \frac{A}{2 T} \lambda \hat{\rho}(0)+{\mathcal
O}\left(\frac{1}{T^2}\right) + {\mathcal O}((\lambda T)^{-\infty}).
\end{multline}
Therefore, one may rewrite $\kappa(\lambda, T)$ as
\begin{multline} \label{prelim2}
\kappa(\lambda,T) = \frac{1}{T} \int_{-\infty}^{\infty}
\hat{\rho}(t/T) \psi(t) \cos(\lambda t) \Tr \,e(t) dt  +\\
\frac{1}{T} \int_{-\infty}^{\infty} \hat{\rho}(t/T) (1-\psi(t))
\cos(\lambda t) \Tr\, e(t) dt- \frac{A}{2 T} \lambda \hat{\rho}(0) +
{\mathcal O}\left(\frac{1}{T^2}\right)
\end{multline}
Consider now the first term on the right-hand side of
(\ref{prelim2}):
\begin{equation}
\label{firstrhs} \frac{1}{T} \int_{-\infty}^{\infty} \hat{\rho}(t/T)
\psi(t) \cos(\lambda t) \Tr \,e(t) dt
\end{equation}
It follows from the trace formula \cite{DG} that the contributions
to \eqref{firstrhs}  from the non-trivial periods $t=L_{\gamma} \neq
0$ with $0 < L_{\gamma} < T_{0}$ are ${\mathcal O}(1/T)$. At the
same time, the wave trace at $t=0$ has the singularity expansion
\cite{Zelditch}: $\Tr e(t) = -\frac{1}{2\pi} t^{-2}+a_1 t^{-1}+ a_2
+\dots$, where the leading coefficient can be computed by
integrating the principal on-diagonal term of the parametrix for the
wave kernel, see \cite[section 3.1]{JP}.

Taking the contribution  of singularity at zero into account, we
obtain that \eqref{firstrhs} can be represented as
\begin{equation} \label{smalltime}
  \frac{A}{2T} \hat{\rho}(0) \lambda +{\mathcal O}\left(\frac{1}{T}\right).
\end{equation}
Therefore, the leading term in (\ref{firstrhs}) cancels the
$\lambda$-term in (\ref{prelim2}) (cf. \cite[Lemma 3.2.1]{JP} where
a local analogue of such a  result is established). To complete the
proof of the lemma we note that by definition $\chi(t,T) =
\hat{\rho}(t/T) (1-\psi(t))$.  \qed

\subsection{Preliminaries and notations}
We now briefly recall the calculus of small-scale
$\hbar$-pseudodifferential operators \cite{DS,Sj} that will be
needed to carry out the various microlocalizations. In the following
we use the notation $\hbar=\lambda^{-1}$.

Given $0 \leq \delta < \frac{1}{2}$, we say that $a(x,\xi;\hbar) \in
C^{\infty}_{0}(T^{*}X \times (0,\hbar_{0}] )$ is in the symbol class
$S^{m}_{\delta}(T^{*}X)$ provided
$$
|\partial^{\alpha}_{x} \partial^{\beta}_{\xi} a(x,\xi;\hbar)| \leq
C_{\alpha, \beta} \hbar^{-m - \delta ( |\alpha| + |\beta|)}.
$$
Given $a \in S^{m}_{\delta}(T^{*}X)$, one can define the
corresponding $\hbar$-pseudodifferential operator $Op_{\hbar}(a)$
invariantly in terms of the Schwartz kernel:
\begin{equation} \label{pdo}
Op_{\hbar}(a)(x,y) = (2\pi \hbar)^{-n} \int_{T^{*}_{x}X}
e^{-i \exp^{-1}_{x}(y) \cdot \xi/\hbar } \, a(x,\xi;\hbar) \, \zeta (r^{2}(x,y)) d\xi.
\end{equation}
Here, $\exp:T_{x}X \rightarrow X$ is the geodesic exponential map, $\xi \in T_{x}^{*}X, \, r(x,y)$ is geodesic distance between $x$ and $y$ and $\zeta \in C^{\infty}_{0}({\mathbb R})$ is supported in a ball
 $B_{\epsilon}(0)$ and equal to $1$ in $B_{\epsilon/2}(0)$ with $\epsilon >0$ sufficiently small
(one can take here $\epsilon < {\rm inj}(X,g)$. Such operators
form a calculus with $Op_{\hbar}(S^{m}_{\delta}) \circ
Op_{\hbar}(S^{m'}_{\delta}) \subset Op_{\hbar}(S^{m +
m'}_{\delta})$. Calderon-Vaillancourt $L^{2}$-boundedness  and the
$\hbar$-Egorov theorem also hold \cite[section 2]{Sj}. Moreover,
let $x,y \in {\mathbb R}^{n}$ be the coordinates of the points
$x,y \in X$ in some local coordinate system
(here we abuse notation slightly and denote points on the manifold
and their  coordinates by the same letters). By using  the Taylor
expansion $- \exp^{-1}_{x}(y) = x-y + {\mathcal O}(|x-y|^{2}),$ to
make the Kuranishi change of variables $ \xi \mapsto (1 +
{\mathcal O}(x-y)) \xi$ in (\ref{pdo}) and integrating by parts in
$\xi$, one can locally rewrite (\ref{pdo}) in the somewhat more
familiar form $Op_{\hbar}(a)(x,y) = (2\pi \hbar)^{-n}
\int_{{\mathbb R}^{n}} e^{i \langle x-y, \xi \rangle /\hbar } \,
a(x,\xi;\hbar)( 1 + {\mathcal O}(\hbar^{1-2\delta}))  \, \zeta (|x-y|^{2})
d\xi.$ However, it is useful here  to work with the invariantly
defined operators in (\ref{pdo}), and we will do so without
further comment.

 We will also use the following geometric notations. Let $M$ be
the universal cover of $X$. The fundamental group of $X$ is denoted
by $\Gamma$. For $\omega \in \Gamma$, let $L_\omega=L_\gamma$  be
the length of the unique closed geodesic $\gamma$ on $X$ that
corresponds to the conjugacy class $[\omega]\subset \Gamma$. The
lifted bicharacteristic curve in $T^*M$ projecting to $\gamma$ is
denoted by $\tilde{\gamma} \subset T^{*}M$.  The standard cotangent
projection map is denoted by $\pi: T^{*}M \rightarrow M$. We use
analogous notation on $T^{*}X$. We will also denote a periodic
geodesic on $X$ (resp. $T^{*}X$)  and a choice of lift on $M$ (resp.
$T^{*}M$) generally by the same letter when the choice of lift is
uniquely specified. Similarly, functions in $C^{\infty}(T^{*}X)$
will be identified with their lifts to $C^{\infty}(T^{*}M)$.

\subsection{Hadamard parametrix}

One needs to get an asymptotic exponential sum formula for the
leading term in (\ref{split1}) just like in the standard case where
$T(\hbar) \sim 1$.
As we have already indicated, we will need to microlocalize the wave
trace on small-scales with respect to $\hbar = \lambda^{-1}$, so
it's useful to introduce the parameter $\hbar$ in the lifted wave
operator.  We do   this by writing the real part of the wave
operator on $M$ in the form $\tilde{e}(t) = \cos t
\sqrt{\tilde{\Delta}} = \cos t \left ( \frac{ \hbar
\sqrt{\tilde{\Delta}}}{\hbar} \right)$. Here we denote by $\tilde
\Delta$ the Laplacian on $M$ and by $\tilde e(t)$ the real part of
the wave operator on $M$.

Even though introducing $\hbar$ amounts to simply rescaling the
Hadamard parametrix approximations $\tilde{E}_N(t)$ to
$\tilde{e}(t)$ (see \ref{hadamard1.1} below), it's useful to think
of $\tilde{E}_{N}(x,y,t;\hbar);$
as the Schwartz kernels of a family of $\hbar$-Fourier
integral operators (albeit, with trivial dependence on $\hbar$)
where $\hbar \sqrt{\tilde{\Delta}}$ is a classical
$\hbar$-pseudodifferential operator \cite[chapter 7]{DS}. Since $\inj(M) =
\infty$, the Hadamard parametrix approximation
$\tilde{E}_{N}(t,x,y;\hbar)$  is valid for {\em all} $t > 0$ and is
given by the following well-known formula: \cite{Berard}:
\begin{equation} \label{hadamard1.1}
\tilde{E}_{N}(t,x,y;\hbar) = \hbar^{-\frac{1}{2}} \frac{\partial}{\partial t}\,
\int_{0}^{\infty} e^{i\theta [r^{2}(x,y) - t^{2}]/ \hbar} \,
\left( \sum_{k=0}^{N} a_{k}(x,y) \theta_{+}^{-k - 1/2}  \, \hbar^{k} \right)
\, d\theta,
\end{equation}
In the last identity (\ref{hadamard1.1}), $\partial_{t}$ is
understood in the distributional sense. Also, for  $\alpha \in
{\mathbb R}, \, \theta_{\pm}^{\alpha}$ denote the standard
homogeneous distributions \cite[section 3.2]{Ho1} and $a_{k} \in
C^{\infty}(M \times M); \, k=0,1,2,...$. It is also well-known (see
\cite[section 39]{Berard}) that for fixed $ y \in M, \,
\tilde{U}-\tilde{U}_{N} \in C^{N}(M \times {\mathbb R}^{+})$ and
moreover, for any compact $Y \subset M,$ there exists
$C_{1},C_{2}>0$ such that
\begin{equation} \label{hadamard2}
\sup_{y \in Y} | \tilde{e}(t,x,y;\hbar )  -
\tilde{E}_{N}(t,x,y;\hbar)|   \leq C_{2} e^{C_{1} |t|} \, \hbar^{N}.
\end{equation}
In (\ref{hadamard2}) the bound is uniform for $(x,y) \in M\times Y$ with $d(x,y)\le d_0$, where $d_0$ is any constant,
and $C_{j}=C_j(d_0,Y,N);\, j=1,2$.  Then
\cite{Berard,CdV}, one can take
\begin{equation} \label{approxX}
E_{N}(t,x,y;\hbar) = \sum_{\omega \in \Gamma; L_{\omega}  \leq
T(\hbar)} \tilde{E}_{N}(t,x,\omega y;\hbar),
\end{equation}
to be the parametrix approximation to the wave operator $e(t) = \cos
t\sqrt{\Delta}$ on $X$.

Moreover, there exists  an appropriate cutoff function $\eta \in
C^{\infty}_{0}(M;[0,1])$ (see \cite[p. 94]{CdV}, \cite[Lemma
34]{Berard}) with $\diam({\rm supp}\,\,\eta)\leq 2 \diam(X)+1$ so
that
\begin{equation} \label{approxX2}
\Tr e(t) = \int_{M} E_{N}(t,x,x;\hbar) \, \eta(x) \,dvol(x) +
{\mathcal O}(e^{C|t|} \, \hbar^{N}),
\end{equation}
for an appropriate $C>0$.

From now on we put
\begin{equation} \label{pardefn}
E_{N,\omega}(t,x,y;\hbar) := \tilde{E}_{N}(t,x,\omega  y;\hbar).
\end{equation}
and so,
\begin{equation} \label{pardefn2}
E_{N}(t,x,y;\hbar) := \sum_{\omega \in \Gamma; L_{\omega}  \leq
T(\hbar)} E_{N,\omega}(t,x, y;\hbar).
\end{equation}

Then, since $|t| \leq \epsilon |\ln \hbar|,$ the RHS in
(\ref{hadamard2}) is ${\mathcal O}(\hbar^{N( 1- C_{1}\epsilon)})$
and so, by choosing $\epsilon < C_{1}^{-1}$ and $N$ large enough, it
follows from Lemma \ref{reform} that
\begin{eqnarray} \label{split2}
\kappa(\lambda,T) = \frac{1}{2T(\hbar)}  \sum_{id\neq\omega \in
\Gamma} \int_{-\infty}^{\infty} \int_{M} E_{N,\omega}(t,x,x;\hbar)
\,  \eta(x) \, \chi(t;T(\hbar)) \, e^{-it/\hbar} \,
dvol(x) dt \nonumber \\
+ \frac{1}{2T(\hbar)}  \sum_{id\neq\omega \in
\Gamma} \int_{-\infty}^{\infty} \int_{M}
E_{N,\omega}(t,x,x;\hbar) \, \eta(x) \, \chi(t;T(\hbar)) \, e^{it/\hbar} \,
dvol(x) dt \nonumber \\
 + \, {\mathcal O}(T(\hbar)^{-1}).
\end{eqnarray}
\begin{remark}
Note that here and further on the sum over $\omega \in \Gamma$ is
finite due to the presence of the cut-off function
$\chi(t;T(\hbar))$: the summation is taken over elements $\omega \neq id$
that such that $L_\omega \le T(\hbar)$.
\end{remark}
The integral in (\ref{split2}) with the the $e^{\pm it/\hbar}$
appearing  has a  total phase function:
$$ \phi^{\pm}(t,\theta;x,y) = (r^{2}(x,y) - t^{2})\theta  \pm t.$$
Since $\theta \geq 0,$
$$ \partial_{t} \phi^{\pm}(t,\theta;x,y) = -2t \theta  \pm 1,$$
and for any $\alpha \geq 0,$
$$ \partial^{\alpha}_{ t} \chi(t;T(\hbar)) = {\mathcal O}_{\alpha}(1),$$
after first integrating by parts in $\theta$ for the $k\geq1$
terms in (\ref{hadamard1.1}), followed by repeated integration by
parts in $t$ one gets that
\begin{eqnarray} \label{split2.1}
\kappa(\lambda,T) = \frac{1}{2T(\hbar)}  \sum_{id\neq\omega \in
\Gamma} \int_{-\infty}^{\infty} \int_{M}
E_{N,\omega}^{+}(t,x,x;\hbar) \, \,  \chi(t;T(\hbar)) \,
e^{it/\hbar} \,
dvol(x) dt \nonumber \\
+ \frac{1}{2T(\hbar)}  \sum_{id\neq\omega \in \Gamma}
\int_{-\infty}^{\infty} \int_{M} E_{N,\omega}^{-}(t,x,x;\hbar) \, \,
\,  \chi(t;T(\hbar)) \,  e^{-it/\hbar} \,
dvol(x) dt \nonumber \\
 + \, {\mathcal O}(T(\hbar)^{-1}),
\end{eqnarray}
where, $E^{\pm}_{N,\omega}(t,x,y;\hbar)$ is defined to be
\begin{equation} \label{cutoff}
\hbar^{-\frac{1}{2}} \eta(x)
\frac{\partial}{\partial t} \, \int_{0}^{\infty} e^{i\theta
[r^{2}(x,\omega y) - t^{2}]/ \hbar} \, \left( \sum_{k=0}^{N}
a_{k}(x,\omega y) \, \theta_{+}^{-k - 1/2} \, \hbar^{k} \right) \, \psi (-2t\theta
\pm 1) \, d\theta.
\end{equation}
From now on, we denote both the Schwartz kernel and the associated  operator by $E^{\pm}_{N,\omega}$.
We note that the error term on the RHS of  (\ref{split2.1})  is the
sum of the ${\mathcal O}(T(\hbar)^{-1})$-error in (\ref{split2}) and
an ${\mathcal O}(\hbar^{\infty})$-piece that arises from inserting
the cutoff function $\psi(-2t\theta \pm 1)$. To see this, split the
$\theta$-integration  into the regions where $|\theta| \leq 1$ and
$|\theta| \geq 1$. Integrating by parts in $t$ over the first region
contributes an ${\mathcal O}(\hbar^{\infty})$-error in
(\ref{split2.1}). The same thing happens over the region where
$|\theta| \geq 1$ by noting that on this set
$$ |-2t\theta \pm 1 |^{-1} \leq |\theta |^{-1}$$
and that $T(\hbar) = O(|\ln \hbar|)$.

The variables $(t,\theta)$ are roughly-speaking dual and the
presence of the cutoffs $\psi(-2t \theta \pm 1)$ will be exploited
later on when we will need to interchange orders of integration in
the expression for $\kappa(\lambda,T)$. The analysis for  each of
the two integral sums in (\ref{split2.1}) is the same and amounts to
time-reversal. We will denote the first sum in (\ref{split2.1}) by
$\kappa^{+}(\lambda,T)$ and the second by $\kappa^{-}(\lambda,T).$

The next step is to microlocalize
on
shrinking $\hbar$-scales near individual lifts of periodic
geodesics. First, we microlocalize on shrinking scales near
$S^{*}X$.

\subsection{Small-scale microlocalization near $S^{*}X$.}
We will need to exploit the separation of geodesics in phase space
given by the dynamical Lemma \ref{dynlemma}. Since the latter result
applies to unit speed geodesics on $S^{*}X$,  we begin by
$\hbar$-microlocalizing the trace $\kappa(\lambda,T)$ near the
co-sphere bundle $S^{*}X$ on small-scales $\sim \hbar^{\delta}; \, 0
< \delta < 1/2$.   The argument here  is quite standard, but for
completeness, we sketch the proof.

Let $d: T^{*}X \times T^{*}X \rightarrow {\mathbb R}$ be the  distance function defined by:
$$d_{T^*X}((x_1,\xi_1),(x_2,\xi_2))=d_{TX}((x_1,\xi_1^*),(x_2,\xi_2^*)),$$
where $\xi_1^*,\xi_2^*$ are dual vectors to $\xi_1,\xi_2$ and
$d_{TX}$ is the Sasaki distance on $TX$. Lemma~\ref{dynlemma}
extends in a straightforward way to the cotangent bundle endowed
with such a distance function, see Remark \ref{cot}.

Let $\zeta \in C^{\infty}_{0}({\mathbb R})$ be a cutoff function
with $\zeta(x) =1$ for $x$ near $0$ and satisfying ${\rm supp}\,
\zeta \subset [-1,1]$.
 For
$0 < \delta <1/2$ we define
$$
\chi^{\delta}(x,\xi;\hbar)  =
\zeta ( \hbar^{-2\delta} d^{2}( (x,\xi), S^{*}X ) \, ).$$

Then, since $\rho \in S({\mathbb R})$, it follows that, modulo
$O(\lambda^{-\infty})$-errors,  in $\kappa(\lambda,T)$ one can sum
over only the eigenvalues $\lambda_{j}$ satisfying $|\sqrt{\lambda_{j}} -
\lambda | \leq \lambda^{\epsilon}$ for any $\epsilon >0$. Write
$\hbar= \lambda^{-1}$ and from now on, consider only such
eigenfunctions, $\phi_j; j=1,2,...$.  Since $(\hbar \sqrt{\Delta} -
1) \phi_{j} = \lambda_j(\hbar) \phi_{j} ; \,
j=1,2,3,...,$ where $|\lambda_{j}(\hbar) -1| = {\mathcal O}_{\epsilon}(\hbar^{1-\epsilon})$ and the $\hbar$-pseudodifferential operator $P(\hbar) =
\hbar \sqrt{\Delta}-1$ is $\hbar$-microlocally elliptic off $S^{*}X
= \{ (x,\xi) \in T^*X; |\xi|_{g} =1 \}$, by  a parametrix
construction in the calculus $Op_{\hbar}(S^{*}_{\delta}),$  it
follows that

\begin{equation} \label{mass}
 \| Op_{\hbar}( 1 - \chi^{\delta}(x,\xi;\hbar) ) \phi_{j} \|_{L^2} =
{\mathcal O}(\hbar^{\infty}).
\end{equation}
 Since $\sum_{j} \rho (T(\lambda ) (\lambda_{j} + \lambda) ) =
{\mathcal O}(\lambda^{-\infty}),$ it then follows from
(\ref{klamt2})  and (\ref{mass}) that
\begin{eqnarray} \label{microsphere}
\kappa(\lambda,T) &=& \frac{1}{2}\sum_{j} \rho (T(\lambda ) (\lambda_{j} -
\lambda) )  - \frac{A \lambda}{2T} \hat{\rho}(0)  + {\mathcal O}(T^{-2}) \nonumber \\
 &=& \frac{1}{2} \sum_{j} \rho (T(\lambda ) (\lambda_{j} - \lambda)
) \cdot \langle Op_{\hbar}( \chi^{\delta}) \phi_{j}, \phi_{j}
\rangle - \frac{A \lambda}{2T} \hat{\rho}(0)  + {\mathcal O}(T^{-2}).
\end{eqnarray}
 In analogy with the construction of the small-cutoff function $\chi^{\delta} \in C^{\infty}_{0}(T^{*}X)$ we choose $\zeta \in C^{\infty}_{0}({\mathbb R})$ as above with $\zeta(u)=1$
near $u=0$.  We define \linebreak  $\tilde{\chi}^{\delta}(x,\xi;\hbar)
=\zeta (\hbar^{-2\delta} d^{2}((x,\xi), S^*M) ) \cdot \tilde{\chi}_{R}(x)$, where
$\tilde{\chi}_{R} \in C^{\infty}_{0}(M)$ with $\tilde{\chi}_{R} =1$
on ${\rm supp}(\eta)$. Then, from Lemma \ref{reform} and
(\ref{microsphere}) it follows that:
\begin{eqnarray} \label{microsphere2}
\kappa^{\pm}(\lambda,T) =  \frac{1}{2T(\hbar)} \sum_{id \neq \omega
\in \Gamma}
 \int_{-\infty}^{\infty} \int_{M} [ \, Op_{\hbar}
^{*}(\tilde{\chi}^{\delta})\cdot E^{\pm}_{N,\omega}
\cdot  Op_{\hbar} (\tilde{\chi}^{\delta}) \,] (x,x,t;\hbar) \nonumber \\
\times \, \,e^{\pm it/\hbar}  \chi(t; T(\hbar)) \, dvol (x) \, dt +
\, {\mathcal O}(T(\hbar)^{-1}).
\end{eqnarray}

From now on, the cutoff $\tilde{\chi}^{\delta}$ will be included in all
computations and for each $\omega \in \Gamma$, the operators $
E_{N,\omega}^{\pm} : C^{\infty}_{0}(M)\rightarrow
C^{\infty}_{0}(M)$ will be replaced by the microlocalizations
$Op_{\hbar} ^{*}(\tilde{\chi}^{\delta}) \cdot E_{N,\omega}^{\pm} \cdot
Op_{\hbar} (\tilde{\chi}^{\delta}) : C^{\infty}_{0}(M)
\rightarrow C^{\infty}_{0}(M)$. To simplify the writing, we
will continue to denote the latter microlocalized operators simply
by $E^{\pm}_{N,\omega}$.

\subsection{Small-scale microlocalization near periodic geodesics.}

The second part of the small-scale $\hbar$-microlocalization
involves microlocalizing on  shrinking $\hbar$-scales near the lifts
to $T^{*}M$ of  individual periodic geodesics in $T^{*}X$. In light
of the  previous section it suffices to take $\gamma \subset
\varOmega(S^{*}M, \hbar^{\delta})$ for any fixed  $\delta \in
(0,1/2)$.

From   Lemma \ref{dynlemma}, it follows that one can put   {\em
disjoint} tubular neighborhoods $ \varOmega(\gamma, e^{-B
T(\hbar)})$ around all lifts of periodic geodesics $\gamma$ with
periods in the time windows $[t, t +\delta']; T_{0} \leq t \leq
T(\hbar)$. The constant  $B>0$ is  {\em uniform} and depends only on
the curvature pinching condition (\ref{curvature}) .

Fix the lift of a periodic  geodesic $\gamma_{0}$ on $M$ and choose
$\epsilon$ so that $ \delta :=B \epsilon < \frac{1}{2}$. Then,
for $(x,\xi) \in T^{*}M,$ we define the small-scale cutoff functions
\begin{equation} \label{ss1}
\zeta_{\gamma_0}^{\delta} (x,\xi;\hbar) = \zeta \left(
\hbar^{-2\delta} \, d^{2}  ( (x,\xi), \gamma_{0}) ) \right).
\end{equation}
Clearly,
$ \zeta_{\gamma_0}^{\delta}  \in S^{0}_{\delta}(T^{*}M),$ where we
have arranged that $\delta < 1/2$.  Then, for any $\omega \in
\Gamma$, we write
$$E^{\pm}_{N,\omega} = Op_{\hbar}^{*}(\zeta_{\gamma_0}^{\delta})
\cdot E^{\pm}_{N,\omega}
\cdot  Op_{\hbar}(\zeta_{\gamma_0}^{\delta}) +
Op_{\hbar}^{*}(\zeta_{\gamma_0}^{\delta} )  \cdot E^{\pm}_{N,\omega} \cdot
Op_{\hbar}( 1- \zeta_{\gamma_0}^{\delta} ) $$
$$ + Op_{\hbar}^{*} (1- \zeta_{\gamma_0}^{\delta} ) E^{\pm}_{N,\omega}
\cdot  Op_{\hbar}(\zeta_{\gamma_0}^{\delta}  ) +  Op_{\hbar}^{*}( 1-
\zeta_{\gamma_0}^{\delta} ) \cdot E^{\pm}_{N,\omega} \cdot Op_{\hbar}(1-
\zeta_{\gamma_0}^{\delta}  ), $$ and taking traces of both sides, we
split the RHS of (\ref{split2}) into three different integral sums.
Taking into account that $\Tr A = Tr A^{*}$ and $\Tr AB = \Tr BA$,
it follows that
\begin{equation} \label{microlocal1}
\kappa^{\pm}(\lambda;T) = \kappa_{11}^{\pm}(\lambda;T) + 2
\kappa_{12}^{\pm}(\lambda;T) + \kappa_{22}(\lambda;T)+O(T(\hbar)),
\end{equation}
where,
\begin{equation} \label{part1}
\kappa^{\pm}_{11}(\lambda;T) = \frac{1}{2T(\hbar)} \sum_{id \neq
\omega \in \Gamma} \int_{-\infty}^{\infty} \Tr [ E^{\pm}_{N,\omega}  \cdot
Op_{\hbar}(\zeta_{\gamma_0}^{\delta}) \cdot
Op_{\hbar}^{*}(\zeta_{\gamma_0}^{\delta}) ]  \, e^{\pm it/\hbar}
\chi(t;T(\hbar)) \, dt,
\end{equation}
\begin{equation} \label{part2}
\kappa^{\pm}_{12}(\lambda;T) = \frac{1}{2T(\hbar)} \sum_{id \neq
\omega \in \Gamma} \int_{-\infty}^{\infty} \Tr [ E^{\pm}_{N,\omega} \cdot
Op_{\hbar}(\zeta_{\gamma_0}^{\delta}) \cdot
Op_{\hbar}^{*}(1-\zeta_{\gamma_{0}}^{\delta}) ] \, e^{\pm it/\hbar}
 \chi(t;T(\hbar)) \, dt,
\end{equation}
and
\begin{equation} \label{part3}
\kappa_{22}^{\pm}(\lambda;T) = \frac{1}{2T(\hbar)} \sum_{id \neq
\omega \in \Gamma} \int_{-\infty}^{\infty} \Tr [E^{\pm}_{N,\omega} \cdot
Op_{\hbar}(1-\zeta_{\gamma_0}^{\delta}) \cdot
Op_{\hbar}^{*}(1-\zeta_{\gamma_0}^{\delta}) ] \, e^{\pm it/\hbar}
\chi(t;T(\hbar)) \, dt.
\end{equation}
We now estimate each of the integral sums in
(\ref{part1})-(\ref{part3}) separately. Roughly speaking, one should
think of the decomposition in (\ref{part1})-(\ref{part3}) as
follows: (\ref{part1}) gives the microlocal contribution of single
periodic geodesic $\gamma_{0}$ to the trace, (\ref{part2}) consists
of cross-terms which we will show are ${\mathcal O}(\hbar^{\infty})$
and finally, (\ref{part3}) is estimated in the same way as
(\ref{part1}) by successively microlocalizing around all other
periodic geodesics $\gamma \neq \gamma_{0}$ with $L_\gamma \leq
T(\hbar)$. Due to the small-scale microlocalizations and ultimately,
the splitting of the time scale into short time-windows, the
expansions in $\kappa(\lambda, T)$  are no longer classical
polyhomogeneous in $\hbar$. For this reason, it is necessary to give
a somewhat different argument than in the classical $T(\hbar) =
{\mathcal O}(1)$ case (see~\cite{Don}).

\vspace{3mm}

\subsection{Studying $\kappa^{\pm}_{11}(\lambda,T)$ in normal
coordinates} The goal of this section is to prove the formula
\eqref{main1}.
For $(x,y)
\in  \mbox{{\rm supp}} (\eta) \times \mbox{{\rm supp}} (\eta)
\subset M \times M,$ we have
\begin{equation} \label{main0}
E^{\pm}_{N,\omega } \cdot Op_{\hbar}(\zeta_{\gamma_{0}}^{\delta})
\cdot Op_{\hbar}^{*}(\zeta_{\gamma_{0}}^{\delta})(x,y;\hbar) =
\frac{\partial}{\partial t} I_{N,\omega}^{\pm}(x,y,t;\omega,\hbar).
\end{equation}
Here
$$I^{\pm}_{N,\omega}(x,y,t;\omega, \hbar) = (2\pi
\hbar)^{-2}\int\int\int   e^{i[ - t^{2}\theta + r^{2}(x, \omega z)\theta  -
 \exp_{z}^{-1}(y) \cdot \xi
]/\hbar} a_{N}(x,z,\theta;\omega,\hbar) $$
 $$\times \, \eta(x) \,  \tilde{\zeta}_{\gamma_0}^{\delta}
(z,\xi;\hbar) \, \psi(-2t\theta \pm 1) \,  \zeta (r^{2}(z,y)) \, dz
d\xi d\theta,$$ where  $ a_N(x,z,\theta;\omega,\hbar) =
\hbar^{-\frac{1}{2}}  \theta^{- \frac{1}{2}}_{+}\,\sum_{k=0}^{N}
a_{k}(x,\omega z) \hbar^{k} \theta_{+}^{-k} $. By symbolic calculus
for $\hbar$-pseudodifferential operators,
$\tilde{\zeta}_{\gamma_0}^{\delta}  \in S^{0}_{\delta}$  satisfies
(\cite[p. 78]{DS} )
\begin{equation} \label{symbolcalculus}
\tilde{\zeta}_{\gamma_0}^{\delta} (z,\xi;\hbar) = |\zeta_{\gamma_0}^{\delta}(z,\xi)|^{2} +
 \hbar^{1-2\delta} \zeta_{\gamma_0}^{\delta,-1}(z,\xi;\hbar)
\end{equation}
with $|\zeta_{\gamma_0}^{\delta}|^{2}, \zeta_{\gamma_0}^{\delta,-1} \in S^{0}_{\delta}.$  Since ${\rm inj}(M,g) = \infty$, to simplify the writing somewhat, we put here   $\zeta (r^{2}(z,y)) = 1$ in the $\hbar$-pseudodifferential cutoffs (see \eqref{pdo}).


Since $\omega \in \Gamma$ acts by
by isometries on $(M,g)$, it follows that
$ \exp_{\omega z}^{-1}(\omega y) \cdot (d\omega^{-1})^t \xi   = \exp_{z}^{-1}(y) \cdot \xi$. and so  the expression for $I^{\pm}_{N,\omega}(x,y,t;\omega,\hbar)$ above can be rewritten as
$$(2\pi
\hbar)^{-2}\int\int\int   e^{i[ - t^{2}\theta + r^{2}(x, \omega z)\theta -
 \exp_{\omega z}^{-1}(\omega y) \cdot  (d\omega^{-1})^{t} \xi
]/\hbar} a_{N}(x,z,\theta;\omega,\hbar) $$
 $$\times \, \eta(x) \,  \tilde{\zeta}_{\gamma_0}^{\delta}
(z,\xi;\hbar) \, \psi(-2t\theta \pm 1) \, dz d\xi d\theta,$$

Changing the  variables $\omega  z \mapsto z, \, (d\omega^{-1})^{t}
\xi \mapsto \xi$ and using that $\omega \in \Gamma$ acts by
isometries, 
one gets
\begin{equation}
\label{main001}
 I^{\pm}_{N,\omega}(x,y,t;\omega,\hbar) =  (2\pi
\hbar)^{-2}\int\int\int   e^{i[ - t^{2}\theta + r^{2}(x, z)\theta -
 \exp_{ z}^{-1}(\omega y) \cdot   \xi
]/\hbar} a_{N}(x, z,\theta;\hbar) $$
 $$\times \, \eta(x) \,  \tilde{\zeta}_{\gamma_0}^{\delta}
(\omega^{-1} z, (d\omega)^{t} \xi;\hbar) \, \psi(-2t\theta \pm 1) \,
dz d\xi d\theta.
\end{equation}

Here, $\omega$ has been scaled out of the amplitude $a_{N}$ so that $a_{N}(x,z,\theta;\hbar) =\hbar^{-\frac{1}{2} } \theta_{+}^{-\frac{1}{2} } \cdot \sum_{k=0}^{N} a_{k}(x,z) \hbar^{k}
\theta_{+}^{-k}.$

Let us now fix a global normal coordinate system centered at $x\in
M$. With some abuse of notation we identify a point $z \in M$ and a
vector of its coordinates in this system: $z=(z_1,z_2) \in T_x M
\cong {\Bbb R}^2$. For instance, in the following $z-\omega y$ means
a vector $(z_1-(\omega y)_1, z_2-(\omega y)_2)$.

Writing the Taylor expansion $$ - \exp_{ z}^{-1}( \omega y) \cdot
\xi = \linebreak  (z- \omega y) \cdot ( 1 + {\mathcal O}(z- \omega
y)) \xi$$ and  making the corresponding Kuranishi change of
variables $\xi \mapsto \xi( 1+ {\mathcal O}(z- \omega y))$  in the
above integral, we  get a coordinate expression for \eqref{main001}:

\begin{eqnarray} \label{main00}
I^{\pm}_{N,\omega}(x,y,t;\omega, \hbar) = (2\pi
\hbar)^{-2}\int\int\int   e^{i[ - t^{2}\theta + r^{2}(x,z)\theta +
(z-\omega y)\xi
]/\hbar} a_{N}'(x,z,\theta;\omega,\hbar) \nonumber \\
 \times \, \eta(x) \,  \tilde{\zeta}_{\gamma_0}^{\delta}
( \omega^{-1} z, (d\omega)^{t} \xi;\hbar) \, \psi(-2t\theta \pm 1) \,  dz d\xi d\theta,
\end{eqnarray}
where $a_{N}' (x,z,\theta;\omega,\hbar)  = a_{N} (x,z,\theta;\hbar) \cdot ( 1 + {\mathcal O}(z-\omega y)).$
Here $\omega$ (again, with some abuse of notation) is understood as
the transformation ${\Bbb R}^2 \to {\Bbb R}^2$, mapping the
coordinates of $z \in M$ to the coordinates of $\omega z \in M$ in
the normal coordinate system centered at $x\in M$; $d\omega$ is
understood as the Jacobian of this mapping. Note that in formula
\eqref{symbolbound} and Lemma \ref{deriv} below, $\omega$ and
$d\omega$ are also understood in this sense.

Recall from $(\ref{symbolcalculus}) $ that  one can write $\tilde{\zeta}_{\gamma_0}^{\delta} = |\zeta_{\gamma_0}^{\delta}|^{2} + \hbar^{1-2\delta}  \zeta_{\gamma_0}^{\delta,-1}$,  
 where the second term $\zeta_{\gamma_0}^{\delta,-1} \in S^{0}_{\delta}$  satisfies estimates of the form 
\begin{equation} \label{symbolbound}
\partial_{z}^{\alpha} \partial_{\xi}^{\beta} \zeta_{\gamma_0}^{\delta,-1}(\omega^{-1}z, (d\omega)^{t} \xi; \hbar )=
O_{\alpha,\beta}\left(e^{C L_{\omega}} \hbar^{-\delta (|\alpha| +
|\beta|)}\right).
\end{equation}
Here,  $C$  denotes possibly different positive constants not
depending on $\omega$. The constants $C$ may depend on $\alpha$ and
$\beta$,  but this dependence can be ignored since we need to take
into account only a finite number of derivatives: $|\alpha+\beta|\le
2(m+1)$, where $m$ is given by \eqref{mlimit}.

 To prove
\eqref{symbolbound} we first note that since
$\zeta_{\gamma_0}^{\delta,-1} \in S^{0}_{\delta}$, by the chain
rule, differentiating the symbol gives the negative powers of
$\hbar$. The estimate then follows from Lemma \ref{deriv} below,
since  $\omega^{-1} z \in \pi ({\rm supp}\,
\zeta_{\gamma_0}^{\delta,-1})$ and therefore $r(x,\omega^{-1}z)< C$.

\begin{lemma}
\label{deriv} Let $M$ and $\omega$ be as above. Fix a normal
coordinate system centered at $x \in M$,  and let $z \in M$ be such
that $r(x,\omega^{-1}z)< C$. Let $z=(z_1,z_2)$, $\omega^{-1} z =
((\omega^{-1} z)_1, (\omega^{-1} z)_2)$ be, respectively,
coordinates of $z$ and $\omega^{-1} z$ in this system, and let
$\xi=(\xi_1,\xi_2)\in T_z^*M$ and $d\omega^t \xi = ((d\omega^t
\xi)_1, (d\omega^t \xi)_2)\in T_{\omega^{-1} z}^* M$ be the
corresponding covectors. Then $\partial_z^{\alpha}(\omega^{-1}
z)=O_\alpha(e^{C L_\omega})$ and $\partial_\xi^{\beta} (d\omega^t
\xi)=O_\beta(e^{C L_\omega})$.
\end{lemma}
\noindent {\bf Proof.} Since $\omega$ is an isometry, we can
identify the tangent spaces $T_xM \cong {\Bbb R}^2$ and $T_{\omega
x} M \cong {\Bbb R}^2$ using $d\omega$,  so that for any point $y
\in M$, the coordinates of the point $y$ in a normal coordinate
system centered at $x$ coincide with the coordinates of $\omega y$
in the normal coordinate system centered at $\omega x$. Taking this
identification into account we obtain
\begin{equation}
\label{qqq} ((\omega^{-1} z)_1, (\omega^{-1} z)_2)=\exp_{\omega
x}^{-1} \circ \exp_x (z_1,z_2).
\end{equation}
As was shown in \cite[Appendix, Propositions 1 and 3]{Berard}, the
Jacobian of the exponential map is bounded away from zero and the
derivatives of the Jacobi fields have at most exponential growth at
infinity. At the same time, the derivatives of the exponential map
could be expressed in terms of the derivatives of the Jacobi fields
\cite[p. 103]{Chavel}, and hence also have at most exponential
growth at infinity. Therefore, the derivatives of the  map
$\exp_{\omega x}^{-1}$ in  \eqref{qqq} are ${\mathcal O}(1)$ since
by assumption of the lemma $r(z,\omega x) < C$. The derivatives of
the map $\exp_x$ in \eqref{qqq} grow exponentially in $r(x,z)\leq
r(x,\omega x)+r(\omega x, z) \leq r(x,\omega x) + C$. At the same
time, $r(x,\omega x) \le L_\omega + C$ by triangle inequality. This
completes the proof of the first estimate in Lemma \ref{deriv}.

Consider now the $\xi$-derivatives and let $\beta=(\beta_1,\beta_2)$. Note that
$\partial_\xi^{\beta} (d\omega^t \xi)=0$ if $\beta_1 \ge 2$ or
$\beta_2 \ge 2$. The first derivatives of $d\omega^t \xi$ pull out
components of $d\omega^{t}$, which is the transpose matrix of the
differential $d\omega$.  Let us recall in what sense $d\omega$ is
understood here: it is the Jacobian of the map $\omega: {\Bbb R}^2
\to {\Bbb R}^2$ given by \eqref{qqq}. But we have already proved
above  that the derivatives of this map (in particular, the first
derivatives that are components of the Jacobian) grow at most
exponentially in $L_\omega$. This completes the proof of Lemma
\ref{deriv}. \qed

\smallskip

Let us go back to the formula (\ref{main00}).  By integration by
parts in $\xi$ in (\ref{main00}), modulo ${\mathcal
O}(\hbar^{\infty})$-errors, one can assume that $r(z, \omega y) =
r(\omega^{-1}z,y) \leq \frac{1}{C}$.  By carrying out the same
argument as for the leading term $|\zeta_{\gamma_0}^{\delta}|^{2}$
below and using the derivative bounds \eqref{symbolbound}, it
follows that by taking $\delta
>0$ sufficiently small,  the contribution to
$\kappa_{11}(\lambda,T)$ of the remainder term $ \hbar^{1-2\delta}
\zeta_{\gamma_0}^{\delta,-1}(\omega^{-1}z, (d\omega)^{t} \xi;\hbar)$
is ${\mathcal O}(\hbar^{C(\delta)})$ for some $C(\delta) >0$.  From
now on, with a slight abuse of notation, we ignore this remainder
and rewrite $|\zeta_{\gamma_0}^{\delta}|^{2}$ as
$\tilde{\zeta}_{\gamma_0}^{\delta}$. So, from (\ref{main00}), we
need to study the integral
 \begin{eqnarray} \label{main1}
I^{\pm}_{N,\omega}(x,y,t;\omega, \hbar) = (2\pi
\hbar)^{-2}\int\int\int   e^{i[ - t^{2}\theta + r^{2}(x,z)\theta +
(z-\omega y)\xi
]/\hbar} a_{N}'(x,z,\theta;\omega,\hbar) \nonumber \\
 \times \, \eta(x) \,  \tilde{\zeta}_{\gamma_0}^{\delta}
(  \omega^{-1} z, (d\omega )^{t} \xi;\hbar) \, \psi(-2t\theta \pm 1) \,  dz d\xi d\theta.
\end{eqnarray}

By integration by parts in the $\theta$-variable in (\ref{main1})
it suffices modulo ${\mathcal O}(\hbar^{\infty})$-errors to assume
that for any fixed $\epsilon>0,$
\begin{equation} \label{range}
|r^{2}(x,z) - t^{2}| < \epsilon.
\end{equation}
\subsection{Stationary phase in $(z,\xi)$-variables and an expansion for
$\kappa^{\pm}_{11}(\lambda,T)$} The goal of  this section is to
prove Lemma \ref{main3.9}. We would like to apply  stationary phase
with parameters in the $(z,\xi)$-variables in (\ref{main1}) (see
\cite{Ho1} Theorem 7.7.5). For this, we need the following estimates
for derivatives of the phase function $\phi(z,\xi;x, y, \omega) :=
\theta r^{2}(x, z) + (z- \omega y)\xi $.   First, we have that for
$(\omega^{-1}z, d\omega^{t}\xi;x,y, \omega) \in {\rm
supp}(\tilde{\zeta}^{\delta}_{\gamma_{0}}) \times {\rm supp}(\eta)
\times {\rm supp}(\eta) \times \{ \omega; L_{\omega} \leq T(\hbar)
\},$
\begin{equation} \label{main1.1}
\partial_{z}^{\alpha} \partial_{\xi}^{\beta} \phi(z,\xi;x, y,
\omega) = {\mathcal O}_{\alpha,\beta} ( T(\hbar)^{2} ).
\end{equation}
The estimate (\ref{main1.1}) follows from  the fact that
$T(\hbar)^{-1} \leq \theta \leq T_{0}^{-1}$ on {\rm supp}
$\psi(-2t\theta \pm 1)$ and the following bounds:
$$\partial_{z}^{\alpha} \phi(z,\xi;x,y,\omega) = \partial_{z}^{\alpha}
( \theta r^{2}(x,z)) = \theta \partial_{z}^{\alpha} r^{2}(x,z) = \theta \partial^{\alpha}_{z} (|z|^{2}) $$
$$ = {\mathcal O}_{\alpha} ( \theta \, r^{2}(x,z)) = {\mathcal
O}_{\alpha}(1) T(\hbar)^{2}; \,\,\,\,|\alpha| \geq 0.
$$
Here, we have used (\ref{range}) as well as the fact that
 $r^{2}(x,z) = z_{1}^{2} + z_{2}^{2}$ since the $z_j$'s are  geodesic normal
  coordinates at $x \in M$.   The mixed
$(z,\xi)$-derivatives of $\phi$ are pointwise ${\mathcal O}(1)$.
Finally,
$$
\partial_{\xi} \phi (z,\xi;x,y,\omega) = z- \omega y = {\mathcal O}
(L_{\omega}) = {\mathcal O}( T(\hbar)),
$$
where the last estimate follows from the triangle inequality for the
distance function. Moreover, the ${\mathcal
O}_{\alpha,\beta}(1)$-constants appearing on the RHS in
(\ref{main1.1})  are all uniform in $\omega \in \Gamma$.

The required lower bound on the norm of the gradient says  that for
$(\omega^{-1}z, d\omega^{t}\xi,x,y, \omega) \in {\rm
supp}(\tilde{\zeta}^{\delta}_{\gamma_{0}}) \times {\rm supp}(\eta)
\times {\rm supp}(\eta) \times  \{ \omega; L_{\omega} \leq T(\hbar)
\}, $ there exists another constant $C>0$ (uniform in all
parameters including $\omega$), such that
\begin{equation} \label{main1.2}
| \partial_{z} \phi |^{2} (z,\xi;x,y)  + |\partial_{\xi}\phi|^{2}
(z,\xi;x,y) \geq \frac{1}{C} ( |\xi + \partial_{z}r^{2}(x,z)
\theta|^{2} + |z- \omega y|^{2}),
\end{equation}
Indeed, the lower bound in (\ref{main1.2})  follows by Taylor
expansion around  the critical point $\xi = -\partial_z r^{2}(x,z) \theta, z= \omega y$
and the Hessian lower bound $ ( \partial^{2}_{z,\xi} \phi )  \gg
\, Id$. The latter follows from the fact that   $
\partial_{z}^{2} \phi = \theta \partial_{z}^{2} r^{2}(x,z)) = \theta \partial_{z}^{2} (z_{1}^{2} + z_{2}^{2}) = {\mathcal O}(\theta) =
{\mathcal O}(1)$, $\partial_{z} \partial_{\xi} \phi =1$
and finally, $\partial_{\xi}^{2} \phi = 0$. Given (\ref{main1.1})
and (\ref{main1.2}), by H\"{o}rmander's interpolation proof of
stationary phase (\cite{Ho1} Theorem 7.7.5),
\begin{equation} \label{main2}
\begin{aligned}
I^{\pm}_{N,\omega}(x,y,t;\hbar) &= (2\pi \hbar)^{-1/2} \,
\sum_{k=0}^{N}  \int_{0}^{\infty} e^{i [-t^{2} + r^{2}(x, \omega
y)]\theta/\hbar}\\
&\times b_{N}^{k}(x,y,\theta;\omega, \hbar) \, \psi (-2t\theta \pm
1) \, \theta_{+}^{-1/2 -k} d\theta.
\end{aligned}
\end{equation}
Set
\begin{equation}
\label{mlimit}
 m= \left[ \frac{4}{1-2\delta} \right]
\end{equation}
In (\ref{main2}) we have
\begin{eqnarray} \label{symbol0}
b_{N}^{k}(x,y,\theta;\omega,\hbar) = \hbar^{k}  \sum_{j=0}^{m}
\frac{ \hbar^{j} (D_{z}D_{\xi})^{j}}{(j+1)!} [a_{k}(x,z) \cdot
\tilde{\zeta}_{\gamma_{0}}^{\delta}(\omega^{-1}z, d\omega^{t}
\xi)]|_{\xi =
-\partial_{z}r^{2}(x,z)\theta, z= \omega y} \nonumber \\
+ {\mathcal O} ( \hbar^{k+m-1 } \, T(\hbar)^{2} \, \sup_{|\alpha +
\beta| \leq 2(m+1)} |\partial_{z}^{\alpha}  \partial_{\xi}^{\beta} (
a_{k} \cdot \tilde{\zeta}_{\gamma_{0}}^{\delta})| ).
\end{eqnarray}


Let us now estimate the contribution to
$\kappa_{11}^{\pm}(\lambda,T)$ coming from the error in
(\ref{symbol0}). Since it is well-known that $\partial_{z}^{\alpha}
a_{k} =  {\mathcal O}_{\alpha}(e^{C L_{\omega}})$  (\cite[Appendix,
Lemma 4]{Berard}), we have that
\begin{eqnarray} \label{error1}
\hbar^{m-1}  \sup_{|\alpha + \beta| \leq 2(m+1)}
|\partial_{z}^{\alpha}
\partial_{\xi}^{\beta} ( a_{k} \cdot \tilde{\zeta}_{\gamma_{0}}^{\delta})
| \ll e^{CL_\omega} \hbar^{-2 (m+1)\delta + m -1}) \nonumber \\
\ll e^{C\epsilon |\ln \hbar|} \, \hbar^{(m+1)(1-2\delta) -2 } \ll
\hbar^{(m+1)(1-2\delta) -2 - C \epsilon }\ll \hbar^{2-C\epsilon}
\end{eqnarray}
Consequently, from (\ref{main2}), (\ref{main1}) and (\ref{part1})
it follows that the contribution to $\kappa_{11}^{\pm}(\lambda,T)$
coming from the remainder in (\ref{symbol0}) for each $k=0,...,N,$ is bounded by
\begin{multline} \label{remainder}
\frac{T(\hbar)^{2}}{T(\hbar)} \sum_{T_{0} \leq L_\omega \leq
T(\hbar)} \int_{-\infty}^{\infty} \int_{0}^{\infty} | \psi(-2t\theta
\pm 1)| \chi(t;T(\hbar)) |\partial_{t} e^{\pm it/\hbar}| \hbar^{2- C
\epsilon }\, \hbar^{k}\, \theta_{+}^{-1/2 -k} \, d\theta dt
\\ \ll  \hbar^{k +1- C \epsilon} T(\hbar) \,
 \sum_{L_\omega \leq T(\hbar)} \int_{T_{0}}^{T(\hbar)}
 \int_{1/T(\hbar)}^{1/T_{0}} \theta_{+}^{-1/2 -k} \, d\theta \,dt
 \\
\ll \hbar^{k +1-C \epsilon} \, T(\hbar) \,e^{C|T(\hbar)|}
 \int_{T_{0}}^{T(\hbar)} |T(\hbar)|^{-1/2 +k} \, dt
\ll \,\hbar^{k+1- C \epsilon} \, |T(\hbar)|^{3/2 +k}.
\end{multline}
In the second to last estimate in (\ref{remainder}) we have used the
exponential  bounds for the growth of the number of periodic geodesics  \cite{MS}.

Choosing $\epsilon>0$ small enough and taking into account that
$T(\hbar) \leq \epsilon |\ln \hbar|$, it follows that the error
terms for $k \geq 0$ in (\ref{symbol0}) are ${\mathcal
O}(\hbar^{\alpha_{1}})$ for some $\alpha_{1} >0$.
We are only interested here in working up to such errors. So, it is
enough to consider only the principal term in (\ref{symbol0}).
Next, note that for each $j \geq 1$,
\begin{equation} \label{remainder2}
\hbar^{j} \frac{(D_{z}D_{\xi})^{j}}{(j+1)!} [a_{k}(x,z) \cdot \tilde{\zeta}_{\gamma_{0}}^{\delta}(\omega^{-1}z
, d\omega^{t} \xi)]|_{\xi =-\partial_{z}r^{2}(x,z)\theta, \, z = \omega y} =
{\mathcal O}(\hbar^{j(1-2\delta)} e^{C_{j} L_\omega}),
\end{equation}
and so again after possibly shrinking $\epsilon >0$ in  the
$T(\hbar) \leq \epsilon \ln \hbar$, it suffices modulo ${\mathcal
O}(\hbar^{\alpha_{2}})$ to restrict  the analysis to the $j=0$ case.
By the same token, it suffices to restrict to the $k=0$ case in
(\ref{main2}).

To get $\kappa_{11}^{\pm}(\lambda,T)$, we put $y = x$ in
(\ref{main2}), integrate over $x \in \mbox{{\rm supp}} (\eta)$ and
sum over $\omega \in \Gamma; T_{0} < L_\omega < T(\hbar)$.
Substituting the formula (\ref{main2})  in (\ref{part1}) and taking
into account the estimate in (\ref{remainder2}), for appropriate
$\alpha_{3} >0,$ we get the following
\begin{lemma}
\label{main3.9}
\begin{eqnarray*}
\kappa_{11}^{\pm}(\lambda,T) = \frac{1}{2T(\hbar)} \sum_{id \neq \omega \in \Gamma}
  \int_{0}^{\infty} \int_{M}
\int_{-\infty}^{\infty}   e^{i\theta [-t^{2} + r^{2}(x,\omega
x)]/\hbar}   b_{N}^{0}(x,x,\theta;\omega, \hbar)
(\hbar \theta_{+})^{-\frac{1}{2}}  \, \nonumber \\
\times  \eta(x) \, \psi (-2t\theta\pm1) \frac{\partial}{\partial t}
 \left[  e^{\pm i t/\hbar} \,
 \chi(t;T(\hbar))   \right] \, dt dvol(x)  d\theta  +
 {\mathcal O}(\hbar^{\alpha_{3}}) \nonumber \\
= \frac{1}{2 \hbar T(\hbar)} \sum_{id \neq \omega \in \Gamma}
\int_{0}^{\infty} \int_{M} \int_{-\infty}^{\infty}
e^{i\theta [- t^{2} + r^{2}(x,\omega x) \pm t]/\hbar}
b_{N}^{0}(x,x,\theta;\omega, \hbar)
(\hbar \theta_{+})^{-\frac{1}{2}}  \,  \nonumber \\
\times \eta(x) \, \psi (-2t\theta\pm 1)  \chi(t;T(\hbar))    \, dt dvol(x)  d\theta
+ {\mathcal O}(\hbar^{\alpha_{3}}).
\end{eqnarray*}
\end{lemma}
In the last line in (\ref{main3.9}), we have used that the leading
term comes from applying $\partial_{t}$ to the exponential $e^{\pm
it/\hbar}$ since $\partial_{t} \chi(t;T(\hbar)) = \partial_{t}  [
\hat{\rho}(t/T) \, (1-\psi)(t) ] = {\mathcal O}(1)$.

For a fixed $\omega \in \Gamma$, every term in (\ref{main3.9}) can
be represented as a sum using the expansion \eqref{symbol0} for
$b_N^0$. Let us estimate each term in this sum separately. Given the
estimate in (\ref{remainder2}) the $j=0$ -- term dominates in
\eqref{symbol0} and so, up to ${\mathcal
O}(\hbar^{\alpha_{4}})$-errors, it suffices to estimate:
\begin{equation} \label{symbol4.1}
\int \int \int e^{i \theta [- t^{2} + r^{2}(x,\omega x)]/\hbar \pm it /\hbar}
 \,\tilde{\zeta}_{\gamma_0}^{\delta}(x,
- d\omega^{t} \cdot \partial_{z}r^{2}(x,z)|_{z=\omega x} \, \theta)
\end{equation}
$$\times  \eta(x) \, a_{0}(x,\omega x) \,
(\hbar \theta_{+})^{-\frac{1}{2}}  \, \chi(t;T(\hbar))\,
\psi(-2t\theta \pm 1)
\,   dt \,dvol(x) \, d\theta. $$

\subsection{\bf Splitting into short time-windows.}
The goal of this section is to prove  Lemma \ref{final0}. We split
up the iterated time integral in (\ref{symbol4.1}) into a sum over
'time-windows' of size $\frac{\delta'}{T(\hbar)}$. The reason for
this is that Lemma \ref{dynlemma} only controls the splitting of
geodesics in a time interval of size $\frac{\delta'}{T(\hbar)}$
(note that the original time variable $t$ has already been rescaled
to $\frac{t}{T(\hbar)}$ at this point). Consider a covering of
$[T_{0}, T(\hbar)]$ by open intervals $I_j$ of length $\delta'$,
$j=0,\dots [\frac{T(\hbar)-T_0}{\delta'}]$,
and let $\eta_j \in C^{\infty}_0 ({\mathbb R})$
be a partition of unity subordinate to this covering. For
$j=0,....,[\frac{T(\hbar)-T_0}{\delta'}]$, we  will need the
following cutoff functions:
\begin{equation} \label{timecut}
\chi_{j}(t;T(\hbar)) = \eta_{j}(T(\hbar) t) \cdot \chi(t;T(\hbar)) \nonumber
\end{equation}
\begin{equation} \label{timecut2}
  = \eta_{j}(T(\hbar) t) \cdot \hat{\rho} ( t/T(\hbar))  \cdot (1-\psi)(t).
\end{equation}

Since we have already inserted the cutoff function $\psi(-2t\theta
\pm 1)$ in  (\ref{symbol4.1}) this allows us to apply Fubini and  do
the $(t,\theta)$-iterated integrals  in (\ref{symbol4.1}) first and
the $x$-integration last. Indeed, we just rewrite the total  phase
in (\ref{symbol4.1}) as a sum:
\begin{equation} \label{phasesplit}
r(x,\omega x) + \Phi^{\pm}(t,\theta;x,\omega),
\end{equation}
where,
$$\Phi^{\pm} (t,\theta;x,\omega) = (\pm t - r(x,\omega x)) -
\theta \, (t^{2} - r^{2}(x,\omega x)).$$ At this point, we would
like to do stationary phase in $(t,\theta)$,  treating $(x,\omega,
\hbar)$ as parameters. Just as in (\ref{main1}) we need to establish
a couple of estimates for derivatives of the phase $\Phi^{\pm}$. The
first estimate
\begin{equation} \label{phasesplit1}
\partial_{\theta}^{\alpha} \partial_{t}^{\beta} \Phi^{\pm}(t,\theta;x,\omega) =
{\mathcal O}_{\alpha,\beta}( r^{2}(x,\omega x))
\end{equation}
with  ${\mathcal O}_{\alpha,\beta}$-constants uniform in all
parameters is immediate.  Since $\partial_{t} \partial_{\theta}
\Phi^{\pm} = -2t, \, \partial_{\theta}^{2} \Phi^{\pm} =0$ and
$\partial_{t}^{2} \Phi^{\pm} = -2 \theta$ one gets the following
lower bound for the $(t,\theta)$- Hessian of $\Phi^{\pm}$:
$$
(\partial_{t,\theta}^{2} \Phi^\pm)^{t} \cdot
(\partial_{t,\theta}^{2} \Phi^\pm)\geq
 \frac{1}{C} \, Id,
$$
and so, by Taylor expansion around the critical points $t = \pm r(x,\omega x)$
and $\theta = \frac{1}{2 r(x,\omega x)}$, one gets the uniform lower bound
\begin{equation} \label{hessian}
|\partial_{t} \Phi^{\pm} |^{2} + |\partial_{\theta} \Phi^{\pm}|^{2} \geq
\frac{1}{C}  \left(  |t \pm r(x,\omega x)|^{2} + |\theta - \frac{1}{2
r(x,\omega x)} |^{2} \right).
\end{equation}

So, by \cite{Ho1} Theorem 7.7.5, it follows that for the expression in
(\ref{symbol4.1}):
\begin{eqnarray} \label{lead}
\int\int\int e^{i \theta [- t^{2} + r^{2}(x,\omega x)]/\hbar \pm it
/\hbar} \,\tilde{\zeta}_{\gamma_0}^{\delta}(x,
- d\omega^{t} \cdot \partial_{z}r^{2}(x,z)|_{z=\omega x} \, \theta
) \cdot \eta(x) \, a_{0}(x,\omega x) \,  \nonumber \\
\times
(\hbar \theta_{+})^{-\frac{1}{2}} \chi(t;T(\hbar))\, \psi(-2t\theta \pm 1)
\,   dt \,d\theta \, dvol(x) \nonumber \\
= 2\pi \hbar^{1/2} \int_{M} e^{ir(x,\omega x)/\hbar}  \eta(x) \,a_{0}(x,\omega x)
\cdot  \tilde{\zeta}_{\gamma_0}^{\delta}(  x, - d\omega^{t} \cdot \partial_{z}r (x,z)|_{z=\omega x}) \, \cdot
  \, r(x,\omega x)^{-\frac{1}{2}} \,  \nonumber \\
 \times \psi(\pm r(x,\omega x)/T) \,  (1-\psi) (\pm r(x,\omega x) ) \, dvol(x) \,
  + {\mathcal O} ( \hbar^{2-2\delta} \, e^{CL_\omega} ).
\end{eqnarray}
In the last line of (\ref{lead}), we have again used the exponential
bounds for the derivatives of  $a_{0}$ \cite{Berard} as well as the
uniform lower bounds for $\Phi^{\pm}$ in (\ref{hessian}). The
contribution of the $ {\mathcal O} ( \hbar^{2-2\delta} \,
e^{CL_\omega})$  error term to $\kappa_{11}(\lambda,T)$ is then
$$
\ll \frac{\hbar^{2(1-\delta)}}{\hbar T(\hbar)} \sum_{\omega; L_{\omega} \leq T(\hbar)}
 e^{CL_\omega}  \ll \hbar^{1-2\delta} e^{C|T(\hbar)|}
|T(\hbar)|^{-2}.
$$
After possibly reducing the size of $\epsilon >0$ further, this term
is then  ${\mathcal O}(\hbar^{\alpha_{5}})$ for some $\alpha_{5}>0$.

We summarize what was shown so far in the following
\begin{lemma}
\label{final0}
\begin{multline*}
\kappa_{11}^{\pm}(\lambda,T)= \frac{\hbar^{-1/2}}{2 T(\hbar)}
\sum_{id \neq \omega \in \Gamma}
\sum_{j=0}^{[\frac{T(\hbar)-T_0}{\delta'}]} \int_{M} e^{ir(x,\omega
x)/\hbar} a_{0}(x,\omega x)  \tilde{\zeta}_{\gamma_0}^{\delta}(x,
-d\omega^{t} \cdot \partial_{z}r (x,z)|_{z = \omega x} )\\ \cdot \chi_{j}( \pm
r(x,\omega x);T(\hbar))  \,
 r(x,\omega x)^{-\frac{1}{2}} \, \eta(x) \, d vol(x)
+ {\mathcal O}(\hbar^{\alpha_{5}}),
\end{multline*}
for some $\alpha_{5}>0$, where $\eta \in C^{\infty}_{0}(M; [0,1])$
is an appropriate cutoff function with $\diam ({\rm supp} \,\,\eta)
\le 2 \diam (X) + 1$.
\end{lemma}

\subsection{Stationary phase in the $x$-variables} The last step involves expanding each term in the
$\omega$-sum in (\ref{final0})  separately. For fixed  $id \neq
\omega \in \Gamma$, one carries out  stationary phase in
(\ref{final0}) in the $x$-variables transverse to the lift of the
periodic geodesic on $X$ given by
\begin{equation} \label{crit1}
\gamma(\omega) = \{ x \in \mbox{{\rm supp}} (\eta) \subset M;
\nabla_{x} f_{\omega}(x) = 0 \,  \},
\end{equation}
where $ f_{\omega}(x) = r(x,\omega x)$ is the displacement function.

Since the dynamical Lemma \ref{dynlemma} only controls separation
 of geodesics in time intervals of size $\frac{\delta'}{T(\hbar)},$
we estimate the summands:
\begin{multline} \label{summands}
\ \frac{\hbar^{-1/2}}{2 T(\hbar)} \sum_{id \neq \omega \in \Gamma} \int_{M}
e^{ir(x,\omega x)/\hbar} a_{0}(x,\omega x) \cdot
 \tilde{\zeta}_{\gamma_0}^{\delta}(x, - d\omega^{t} \cdot \partial_{z}r (x,z)|_{z=\omega x})
\\
\times \,  \chi_{j}( \pm r(x,\omega^{k} x);T(\hbar))\,
 r(x,\omega x)^{-\frac{1}{2}} \,\eta(x) \, dvol(x) ;
\,\, \, \,\,\,\,j=0,...,[\frac{T(\hbar)-T_0}{\delta'}]
\end{multline}
 separately.
For fixed $id \neq \omega \in \Gamma$ and  $j \in \{
0,...,[\frac{T(\hbar)-T_0}{\delta'}]$, we introduce Fermi
coordinates $(u, s ) \in {\mathbb R} \times [0,L]$ in
(\ref{summands}) (see \cite{CdV}) centered on the geodesic segment
$\gamma(\omega)$ (i.e. $u=0$ on $\gamma(\omega)$)  and apply
stationary phase in the $u$-variables just as in the $T(\hbar) \sim
1$ case (see \cite[section 2]{CdV}). In terms of the
$(u,s)$-coordinates on $M$
\begin{equation} \label{jacobian}
dvol(x) = J(u,s;\omega) \, du ds,
\end{equation}
where, $|\partial_{u}^{\alpha} \partial_{s}^{\beta} J(u,s;\omega)| =
{\mathcal O}_{\alpha,\beta} ( e^{C(\alpha,\beta)L_\omega})$ \cite{CdV}.

In the case at hand, it is necessary  to  control the dependence of
the phase and amplitude of (\ref{final0}) on the parameters,
$\omega$.

We now need the following estimates: The first is a (uniform)
Hessian lower bound (\cite[Lemma 4]{CdV}) which says that
\begin{equation} \label{cdv}
\nabla_{u}^{2} f_{\omega}(u,s) \geq C_{0} >0.
\end{equation}
The constant $C_{0}>0$ in (\ref{cdv}) is {\em uniform} in the
parameters $\omega \in \Gamma$ and depends only the curvature
pinching conditions $-K_{1}^2 \leq K \leq -K_{2}^2$. In addition,
one has the upper bounds \cite{Berard, CdV}:
\begin{multline} \label{cdv2}
\partial_{u}^{\alpha} a_{0}(x,\omega x) = {\mathcal O}_{\alpha}(e^{C L_{\omega}}), \,
\partial_{u}^{\alpha} f_{\omega}(u,s) ={\mathcal O}_{\alpha}( e^{C'
L_{\omega}}),\\
\,  \partial_{u}^{\alpha} J(u,s) = {\mathcal O}_{\alpha} (e^{C''L_{\omega}}),
\end{multline}
where the estimates in (\ref{cdv2}) are all uniform in the Fermi
coordinates $(u,s) \in {\mathbb R} \times [0,L]$. In order to apply
stationary phase with parameters, we need to carefully analyze the
critical sets of the phase function $f_{\omega}(u,s)$ in
(\ref{summands}) for all $\omega \in \Gamma$ with $L_{\omega} \ll
|\ln \hbar|$. We do this, by combining the estimates (\ref{cdv}) and
(\ref{cdv2}) with  the dynamical Lemma \ref{dynlemma}.

\subsection{Application of Lemma \ref{dynlemma}.}
As above, let  $\gamma(\omega)$ be the unique lifted geodesic
invariant under the action of $\omega \in \Gamma$.  As is well-known
\cite[Proposition 4.2]{BO}, the displacement function
$f_{\omega}(x)$ on a negatively curved surface is strictly convex,
except on $\gamma(\omega)$ where it is constant: $f_{\omega}(x) =
L_\omega$ when $x \in \gamma(\omega)$. The geodesic $\gamma(\omega)$
is also the critical set of $f_{\omega}(x)$, see \eqref{crit1}.

Fix $j=0,...,[\frac{T(\hbar)-T_0}{\delta'}]$ and $\omega \in \Gamma$
and consider the corresponding  summand in (\ref{summands}) given by
\begin{equation} \label{summands2}
\frac{\hbar^{-1/2}}{2 T(\hbar)} \int_{M} e^{ir(x,\omega x)/\hbar}
a_{0}(x,\omega x) \cdot \tilde{\zeta}_{\gamma_0}^{\delta}(x,
- d\omega^{t} \cdot \partial_{z}r (x,z)|_{z=\omega x})  \,  \chi_{j}( \pm r(x,\omega^{k}
x);T(\hbar))\,
\end{equation}
$$\times r(x,\omega x)^{-\frac{1}{2}} \,\eta(x) \, dvol(x).$$

We now show that, the integral in  (\ref{summands}) is ${\mathcal
O}(\hbar^{\infty})$ unless $\omega = \omega_{0}$, where the latter
group element fixes the lifted geodesic $\gamma_{0}(\omega_{0})$. As
before, $\pi: T^{*}M \rightarrow M$ denotes the standard cotangent
projection map and $\tilde{\gamma} \subset T^{*}M$ will denote the
bicharacteristic curve with $\pi(\tilde{\gamma}) = \gamma$.

To prove the above claim, we first note that on  {\rm supp} $\tilde
\zeta^{\delta}_{\gamma_{0}}(x, - d\omega^{t} \partial_{z}r (x,z)|_{z=\omega x})$,
\begin{equation} \label{bounds2}
d(x, \pi(\tilde{\gamma_{0}})) = d(x, \gamma_{0})= {\mathcal O}(\hbar^{\delta}).
\end{equation}
Here, following our convention we write $\gamma_{0}$ for
$\gamma_{0}(\omega_{0})$.  Since  $d\omega_{0}^{t} \cdot \partial_{z} r(x,z)|_{z=\omega_{0}x}$
is the cotangent vector to $\tilde{\gamma_{0}}$ at $ x \in
\pi(\tilde{\gamma_0})$, it follows from the small-scale
microlocalization in the $\xi$-variables   that  for any $x \in \gamma_{0},$
we have that
\begin{equation} \label{bounds3-}
d\omega^{t} \cdot \partial_{z} r(x,z)|_{z=\omega x} =  d\omega_{0}^{t} \cdot \partial_{z} r(x,z)|_{z=\omega_{0}x} +
{\mathcal O}(\hbar^{\delta}).
\end{equation}
From the exponential upper  bounds for the derivatives of
$f_{\omega}$ in (\ref{cdv2}) and by a   Taylor expansion around
$\gamma_{0}$, it follows that for any $x \in M$ with
$d(x,\gamma_{0}) = {\mathcal O}(\hbar^{\delta}),$
\begin{equation} \label{bounds3}
d\omega^{t} \cdot \partial_{z} r(x,z)|_{z=\omega x} = d\omega_{0}^{t} \cdot \partial_{z} r(x,z)|_{z=\omega_{0}x} +
{\mathcal O}(\hbar^{\delta_{0}}),
\end{equation}
for some  $ 0 < \delta_{0} < \delta.$

On the other hand, by taking $\epsilon >0$ small enough, and again
using the upper bounds in (\ref{cdv2}), it follows by an
integration by parts in (\ref{final0}) in the transversal $u$-variable, that
modulo ${\mathcal O}(\hbar^{\infty})$-errors, one can cut off the
integration to values of $x$ satisfying:
\begin{equation} \label{crit2}
\partial_{x} f_{\omega}(x) = {\mathcal O}(\hbar^{\delta}),
\end{equation}
for any $0 \leq \delta < 1/2.$  Then using (\ref{cdv2}) again
together with the uniform Hessian  lower bound in (\ref{cdv}), it
follows by a Taylor expansion argument that
\begin{equation} \label{closebase}
d( \gamma_{0}( \omega_{0} ) , \gamma( \omega  ) ) =
{\mathcal O}(\hbar^{\delta_{1}}).
\end{equation}
Here, $0< \delta_{1} < \delta_{0}$ is yet another, possibly smaller constant.
But then,  in view of the estimate (\ref{bounds3}) we get that
\begin{equation} \label{close}
d( \tilde{ \gamma_{0} } ( \omega_{0}   ) ,
\tilde{ \gamma }( \omega  )  ) = {\mathcal O}(\hbar^{\delta_{1}}).
\end{equation}
The presence of the time cutoff $\chi_{j}$ in (\ref{summands})
ensures that there  is only one $j$-summand  in (\ref{summands})
that contributes in a non-negligible way to $\kappa_{11}^{\pm}$,
namely, the interval $I_j$
containing $L_{\gamma_{0}}$.

Finally, by  possibly shrinking $\epsilon >0$ further, it  follows
from (\ref{close})  and  the dynamical Lemma \ref{dynlemma} that, up
to ${\mathcal O}(\hbar^{\infty})$-error in $\kappa_{11}^{\pm}
(\lambda,T)$, $\omega = \omega_{0}.$

One can repeat the above argument for each $j \in
\{0,...,[\frac{T(\hbar)-T_0}{\delta'}]$\} and $id \neq \omega \in
\Gamma$ separately in (\ref{summands2}) and sum up. By taking the
exponential bounds in (\ref{cdv2}) into account and possibly further
decreasing the size of $\epsilon$ in $T(\hbar) = \epsilon |\ln
\hbar|$ if necessary, it follows by applying stationary phase  in
$u$ (see \cite{CdV,Don, Sunada}) that for appropriate
$\alpha_{6}>0$,
\begin{equation} \label{leadestimate}
\kappa_{11}^{\pm}(\lambda, T)  = \frac{e^{\pm i
L_{\gamma_{0}}/\hbar}}{ 2 T(\hbar)} \cdot  L^{\sharp}_{\gamma_{0}}
\cdot \chi \left(L_{\gamma_0},T(\hbar)\right) \, \cdot |\det (I -
P_{\gamma_{0}})|^{-1/2} + E(\hbar),
\end{equation}
where,
\begin{equation}
\label{err1} E(\hbar)={\mathcal O}(\hbar^{\alpha_6})+ {\mathcal O}
\left( \frac{ e^{C T(\hbar)} \hbar^{\infty} }{ T(\hbar)} \right) =
{\mathcal O}(\hbar^{\alpha_{6}}).
\end{equation}
The first term on the RHS of  \eqref{err1} gives the remainder
produced by the stationary phase method for $\omega=\omega_0$. The
second term on the RHS  of  \eqref{err1} follows from the bound $\#
\{ \gamma; L_\gamma \leq T \} = {\mathcal O}(e^{C T(\hbar)}
T(\hbar)^{-1}),$ and the estimate ${\mathcal O}(\hbar^{\infty})$ for
each summand in \eqref{summands} corresponding to $\omega \neq
\omega_0$.

\vspace{3mm} \noindent{\bf Estimate for
$\kappa^{\pm}_{12}(\lambda,T)$:} Lemma \ref{dynlemma}  implies that
there are  {\em no} periodic geodesics in {\rm supp}   $( \,
\zeta_{\gamma_{0}}^{\delta} \cdot (1-\zeta_{\gamma_{0}}^{\delta}  ) \,
)$. So, by repeated  integration by parts in the $u$-variable in
(the analogue of) (\ref{summands2}), it follows that
\begin{equation} \label{errorestimate}
\kappa_{12}(\lambda,T) = {\mathcal O}\left( \frac{ e^{C T(\hbar)}
\hbar^{\infty} }{ T(\hbar)} \right)  = {\mathcal O}(\hbar^{\infty}),
\end{equation}
when $T(\hbar) = \epsilon |\ln \hbar|$.

\vspace{3mm} \noindent{\bf Estimate for
$\kappa^{\pm}_{22}(\lambda,T)$:} Here, repeat the  microlocalization
as in the estimate for $\kappa_{11}(\lambda,T)$ near each periodic
geodesic $\gamma \neq \gamma_{0}$ with $ 0 < L_{\gamma} \leq
\epsilon |\ln \hbar|$ separately  and sum up. By the same argument
as for $\kappa_{12}(\lambda,T)$, all cross terms give ${\mathcal
O}(\hbar^{\infty})$ contributions to $\kappa(\lambda,T)$. Also, we
note that since the remainder term is ${\mathcal
O}(\hbar^{\alpha_{6}})$ for some $\alpha_{6} >0$ in
(\ref{leadestimate}), after summing over  all $\omega \neq
\omega_{0}$ in $\kappa_{22}^{\pm}$, it follows that, after possibly
shrinking $\epsilon >0$ further, the remainder in the latter is
$$
{\mathcal O}(e^{CT(\hbar)} T(\hbar)^{-1} \hbar^{\alpha_{6}}) = {\mathcal O}(\hbar^{\alpha_7}),
$$
for some $\alpha_{7} >0$. This finishes the proof of Proposition
\ref{k2} since the ${\mathcal O}(\hbar^{\alpha_{7}})$-error is
absorbed in the ${\mathcal O}(T^{-1})$-error in (\ref{split1}) for
$\kappa(\lambda,T)$. \qed

\section{Proof of Theorem \ref{main}}\label{sec:mainproof2}
\subsection{Application of thermodynamic formalism}
The strategy of the proof of Theorem \ref{main} is to get a
contradiction with Lemma \ref{above:weak1} that gives an {\em upper
bound} for the quantity $\kappa(\lambda,T)$ defined in
\eqref{klamt2}, \eqref{klamt3}; that bound holds for {\em any}
$\lambda,T$.  We use the formula \eqref{qq} and the estimate
\eqref{grow1} below to obtain a {\em lower bound} for
$\kappa(\lambda,T)$ which holds for an infinite sequence of pairs
$(\lambda,T)$ such that $T\sim \frac{1}{h} \ln\ln\lambda$, where $h$
is  the topological entropy of the geodesic flow on $X$.

Consider the main term in the formula \eqref{qq}.  This sum is a
trigonometric polynomial in $\lambda$, while the number of terms and
the coefficients depend on $T$.

We would like to choose $\lambda$ and $T$ in such a way that the
value of the polynomial is large. We first study the rate of growth
for the sum of coefficients of this trigonometric polynomial as
$T\to\infty$ not taking into account the oscillatory nature of the
terms. Let
\begin{equation}\label{dyn:sum}
S(T)=\sum_{L_\gamma \in \Lsp, L_\gamma\le T}
\frac{L_\gamma}{\sqrt{|\det(I-\cP_\gamma)|}};
\end{equation}

It turns out that the asymptotic rate of growth of $S(T)$ can be
determined using results from the theory of thermodynamic formalism
for Anosov flows:
\begin{prop}\label{sumgrow}
There exists a constant $C_0>0$ such that
\begin{equation}\label{grow1}
S(T)= C_0 e^{P\left(-\frac{\cH}{2}\right)\cdot T}(1+o(1))
\end{equation}
as $T\to\infty$, where $P$ is the topological pressure
\eqref{pressure} and $\cH$ is the Sinai-Ruelle-Bowen potential
\eqref{SRB}.
\end{prop}
It was shown in \cite{JP} that $P\left(-\frac{\cH}{2}\right)\geq
K_2/2$, hence $S(T)$ grows exponentially in~$T$.

\noindent{\bf Proof of Proposition \ref{sumgrow}.} 
Let $\xi_\gamma \in T_{\gamma(0)}M$ be the tangent vector to
$\gamma$. The Poincar\'e map $\cP_\gamma$ preserves the unstable
subspace $E_{\xi_\gamma}^u$ and the stable subspace
$E_{\xi_\gamma}^s$, both of dimension one. The map $\cP_\gamma$ is
symplectic and has determinant equal to one. The eigenvalues
$\mu,1/\mu$ of $\cP_\gamma$, corresponding, respectively, to the
unstable and the stable directions, satisfy:
\begin{equation}\label{eigunstable}
\ln |\mu|  \in[K_2L_\gamma,K_1L_\gamma].
\end{equation}

It follows from \eqref{eigunstable} that
$$
|\det(I-\cP_\gamma)|=(|\mu|-1)(1-\frac{1}{|\mu|}) = |\mu|
\left(1+O(e^{-K_2L_\gamma})\right),
$$
On the other hand, by definition of $\cH$
$$
|\mu|=\exp\left[\int_0^{L_\gamma} \cH(G^s \xi_\gamma)ds\right],
$$
so we have
\begin{equation}\label{normeq}
|\det(I-\cP_\gamma)|=\exp\left[\int_0^{L_\gamma}\cH(G^s\xi_\gamma)ds\right]
\left(1+O(e^{-K_2L_\gamma})\right).
\end{equation}

Split the sum $S(T)$ into two parts:
$$
S(T)=\sum_{L_\gamma\le T/2} +\sum_{T/2< L_\gamma\le
T}=S_1(T)+S_2(T).
$$

It follows from \eqref{dyn:sum} and \eqref{normeq} that
\begin{equation}\label{dyn:sum2}
S_2(T)=\sum_{L_\gamma \in \Lsp, T/2\le L_\gamma\le T} L_\gamma
\exp\left[-\frac{1}{2}\int_0^{L_\gamma}\cH(G^s\xi_\gamma)ds\right]
\left(1+O(e^{-K_2T/2})\right),
\end{equation}
and that
$$
S_1(T)=O\left(\sum_{L_\gamma \in \Lsp, L_\gamma\le T/2} L_\gamma
\exp\left[-\frac{1}{2}\int_0^{L_\gamma}\cH(G^s\xi_\gamma)ds\right]\right).
$$

It follows easily from the results of \cite{Par}, \cite[(7.1)]{PP},
\cite[p. 109]{MS} that
$$
S_2(T)=\frac{e^{P\left(-\frac{\cH}{2}\right)\cdot T}}{P(-\cH/2)}\;
(1+o(1))
$$
and that
$$
S_1(T)=O\left(e^{P\left(-\frac{\cH}{2}\right)\cdot
\frac{T}{2}}\right).
$$
This finishes the proof of the proposition \qed

If follows easily from Proposition \ref{sumgrow} that
\begin{equation}\label{stilde}
\tilde S(t) = \sum_{L_\gamma \in \Lsp, L_\gamma\le T}
\frac{L^\sharp_\gamma}{\sqrt{|\det(I-\cP_\gamma)|}}= C_0
e^{P\left(-\frac{\cH}{2}\right)\cdot T}(1+o(1)).
\end{equation}
Indeed, since each imprimitive geodesic is a multiple of some
geodesic of at least twice smaller length, the contribution of
imprimitive geodesics into $S(T)$ is \newline
$O\left(e^{P\left(-\frac{\cH}{2}\right)\cdot \frac{T}{2}}\right)$
and hence could be neglected (cf. \cite[p. 25]{JP}).

\subsection{Preliminary estimates}
We can now finish the proof of Theorem \ref{main}. Assume for
contradiction that Theorem \ref{main} doesn 't hold. Then
$R(\lambda)$  satisfies
$$
R(\lambda)=O\left((\ln\lambda)^{P(-\cH/2)(1-\varepsilon)/h}\right)
$$
for some $\varepsilon>0$. Let $b=P(-\cH/2)(1-\varepsilon)/h$ be the
exponent $\ln\lambda$ in the previous formula. Then by Lemma
\ref{above:weak1} we have
\begin{equation}\label{ksmall}
\kappa(\lambda,T)=O((\ln\lambda)^b).
\end{equation}

We rewrite \eqref{ksmall} as
\begin{equation}\label{ksmall2}
\ln|\kappa(\lambda,T)|\leq b\ln\ln\lambda+C_1.
\end{equation}
To finish the proof, it suffices to establish a contradiction with
\eqref{ksmall2}.  This will be done using Proposition \ref{k2} and
the estimate \eqref{grow1} for a suitable choice of $\lambda$ and
$T$. In the sequel, we shall let $\lambda,T\to\infty$ while keeping
$T\sim \frac{1}{h}\ln\ln\lambda$.  This ensures that the hypothesis
$T<\epsilon\ln\lambda$ of Proposition \ref{k2} is satisfied.

Denote the main term in the asymptotics of $\kappa(\lambda, T)$ by
$$\Sigma(\lambda,T)=\sum_{L_\gamma
\in \Lsp, L_\gamma\le T} \frac{L^\sharp_\gamma \cos (\lambda
L_\gamma)\,\chi(L_\gamma, T)}{T\sqrt{|\det(I-\cP_\gamma)|}}.
$$ It is a trigonometric polynomial in
$\lambda$. According to \eqref{cutoff0}, without loss of generality
we may assume that
\begin{equation}\label{psi1}
\chi(L_\gamma, T)\geq 1/2,\qquad \forall\,\, L_\gamma \in (T_0,
T(1-\varepsilon/2)].
\end{equation}
The condition $L_\gamma>T_0$ is not essential because it rules out
only a finite number of closed geodesics, and their total
contribution to \eqref{stilde} is $O(1)$.

Next, we would like to choose $\lambda$ so that all the terms
$\cos(\lambda L_\gamma)$ will be (say) $\geq 1/2$ for
$L_\gamma\in\Lsp, T_0< L_\gamma \le T$. Let $\nu(T)$ be the number
of {\em distinct} such $L_\gamma$-s, and let
$L_1,L_2,\ldots,L_{\nu(T)}$ be the corresponding lengths.  It
suffices to choose $\lambda$ so that
\begin{equation}\label{distance1}
{\rm dist}(\lambda L_j,2\pi {\mathbb Z})\ \leq\ 1/2,\qquad 1\leq
j\leq \nu(T).
\end{equation}

Assuming  \eqref{psi1} and \eqref{distance1}, we get for large
enough $T$
$$
\ln|\Sigma(\lambda,T)|\geq C_2+\ln|\tilde S(T)|.
$$
Using the estimate \eqref{stilde} and Proposition \ref{k2} we
conclude that
\begin{equation}\label{klarge1.5}
\ln|\kappa(\lambda,T)| \geq P(-\cH/2)T(1-\varepsilon/2) -\ln T+ C_3.
\end{equation}
This formula will be used to get a contradiction with the upper
bound \eqref{ksmall2} for $|\kappa(\lambda,T)|$.


\subsection{Dirichlet box principle}
We next explain how to choose $\lambda$ so that \eqref{distance1}
would hold.  Let $M_1$ be a large constant whose value will be
specified later.  Then by Dirichlet box principle (\cite{JP}, see
also \cite{PR,RS}) there exists
$$
\lambda\in[M_1,M_1 2^{\nu(T)}]
$$
such that \eqref{distance1} holds.  Hence, for {\em any} choice of
$M_1$ there exists $\lambda$ satisfying
$$
\ln\ln M_1\leq\ln\ln\lambda\leq\ln\ln M_1+\ln \nu(T)+\ln\ln 2.
$$
for which \eqref{distance1} holds.

It follows from results of Margulis \cite{MS} that $\nu(T)=
e^{hT}(1+o(1))/hT$ as $T\to\infty$. Therefore, any $\lambda$
satisfying the previous inequality would also satisfy
\begin{equation}\label{box:ineq}
\ln\ln M_1\leq\ln\ln\lambda\leq \ln\ln(M_1)+hT-\ln(hT)
\end{equation}
We now choose
\begin{equation}\label{M1}
M_1=\exp(\exp(\alpha T)), \quad {\rm where} \,\,\,
\alpha<\frac{h\varepsilon}{2(1-\varepsilon)}.
\end{equation}
Then \eqref{box:ineq} becomes
\begin{equation}\label{box:ineq2}
\alpha T\leq\ln\ln\lambda\leq (h+\alpha)T-\ln(hT).
\end{equation}

The first inequality in \eqref{box:ineq2} ensures that the
hypothesis of Proposition \ref{k2} is satisfied, implying
\eqref{qq}.  By the previous argument, we have shown that
\eqref{box:ineq2} and \eqref{klarge1.5} implies existence of
$\lambda$ such that the following estimate holds:
\begin{equation}\label{klarge2}
\ln|\kappa(\lambda,T)| \geq P(-\cH/2)T(1-\varepsilon/2) -\ln T+ C_3.
\end{equation}

To establish a contradiction with the formula \eqref{ksmall2}, it
suffices to have
$$
P(-\cH/2)T(1-\varepsilon/2) -\ln T+ C_3 >b \ln\ln\lambda+C_1=
\frac{P(-\cH/2)(1-\varepsilon)}{h} \ln\ln\lambda+C_1
$$
or
\begin{equation}\label{contr:ineq}
\ln\ln\lambda\leq \frac{h(1-\varepsilon/2)}{1-\varepsilon}\;
T+\frac{h}{P(-\cH/2)(1-\varepsilon)}(C_4-\ln T)
\end{equation}

If we could show that the inequality \eqref{box:ineq2} {\em implies}
the inequality \eqref{contr:ineq}, we would be done.  Indeed, by
Dirichlet box principle there exists {\em some} $\lambda$ satisfying
\eqref{box:ineq2}, and so \eqref{contr:ineq} holds for {\em that}
value of $\lambda$, establishing a contradiction.  The linear
function of $T$ is the fastest-growing term in the right-hand side
of both inequalities, so it suffices to compare the coefficients of
$T$.  The coefficient in \eqref{box:ineq2} is equal to $h+\alpha$,
while that in \eqref{contr:ineq} is equal to
$h(1-\varepsilon/2)/(1-\varepsilon)$.  Accordingly, it suffices to
have
$$
h+\alpha<h(1-\varepsilon/2)/(1-\varepsilon),
$$
and this is ensured by the choice of $\alpha$ in \eqref{M1}.  This
establishes the desired contradiction and finishes the proof of
Theorem \ref{main}.
\qed

\section{Remainder estimates in higher dimensions}
\subsection{Proof of Theorem \ref{Riesz}} \label{sec:Riesz}
In this section we assume that $X$ is a compact Riemannian manifold
of dimension $n\ge 3$. Recall that the Riesz mean of order $k$ of a
function $f(\lambda)$ is defined by
\begin{equation}
\label{kRiesz}
 \cR_k f(\lambda)=\frac{k}{\lambda} \int_0^{\lambda}
\left(1-\frac{t}{\lambda}\right)^{k-1} f(t) dt, \quad k=1,2,\dots
\end{equation}
As was shown in \cite{Saf2}, on any smooth compact $n$-dimensional
Riemannian manifold
\begin{equation}
\label{Saff}
\cR_2R(\lambda)=C a_1 \lambda^{n-2} + O(\lambda^{n-3}),
\end{equation}
Here $C$ is a non-zero constant depending on the dimension only, and
$a_1=\frac{1}{6}\int_X \tau$, where $\tau$ is the scalar curvature
of $X$. Note that $a_1$ is the first {\it heat invariant} of $X$,
the coefficient in the short time asymptotics of the heat trace:
\begin{equation}
\label{heat} \sum_i e^{-\lambda_i t} \sim \frac{1}{(4\pi)^{n/2}}
\sum_{j=0}^\infty a_j\,\,  t^{j-\frac{n}{2}}.
\end{equation}
Therefore, if $a_1$ does not vanish (which is the assumption of
Theorem \ref{Riesz}), $\cR_2 R(\lambda)>>\lambda^{n-2}$. Combining
\eqref{kRiesz} and \eqref{Saff} we get
$$
\frac{1}{\lambda} \int_0^{\lambda} |R(t)|dt > \frac{1}{\lambda}
\int_0^{\lambda} \left(1-\frac{t}{\lambda}\right) R(t) dt =
\frac{1}{2} \cR_2 R(\lambda)>> \lambda^{n-2}.$$ This completes the
proof of Theorem \ref{Riesz}. \qed
\begin{remark}
As follows from results of \cite{Saf2}, the first Riesz mean of
$R(\lambda)$ satisfies
\begin{equation*}
\frac{1}{\lambda}\int_0^{\lambda} R(t) dt =
O(\lambda^{n-2}).
\end{equation*} At the same time, $R(\lambda)=O(\lambda^{n-1})$.
Putting together these two upper bounds and the lower bound
\eqref{bound3} one gets an idea about  the amount of cancelations
occurring when $R(\lambda)$ is integrated over $[0,\lambda]$.
\end{remark}
\subsection{Oscillatory error term}
\label{oscerr} Following \cite[section 1.2]{JP}  one may  introduce
the {\it oscillatory error term} $R^{osc}(\lambda)$ in Weyl's law:
\begin{equation}
\label{osc} N(\lambda)=\frac{1}{(4\pi)^{n/2}}
\sum_{j=0}^{\left[\frac{n-1}{2}\right]}
\frac{a_j}{\Gamma\left(\frac{n}{2}-j+1\right)}\lambda^{n-2j}+R^{osc}(\lambda),
\end{equation}
where $a_j$ are defined by \eqref{heat}. The expression \eqref{osc}
is {\it not} an asymptotic expansion, however it often appears in
physics literature. Such a representation is quite natural since it
allows to separate the ``mean smooth part'' of the counting function
coming from the singularity of the heat trace at zero, and the
``oscillating part'' produced by the singularities in the wave trace
caused by closed geodesics. We believe that using essentially the
same arguments as in the proof of Theorem \ref{main}, one can show
that the oscillatory error term on an $n$-dimensional compact
negatively curved manifold satisfies:
\begin{equation}
\label{bound1} R^{osc}(\lambda)= \Omega
\left((\ln\lambda)^{\frac{P(-\cH/2)}{h}-\varepsilon}\right) \,\,\,
\forall \, \varepsilon >0.
\end{equation}
In order to prove \eqref{bound1}, one has to extend to dimensions
$n\ge 3$  the dynamical part of the proof of Theorem \ref{main},
which is easy, and to generalize Theorem \ref{tk2}, which requires
some work. In particular, one needs higher-dimensional analogues
of the results of \cite{CdV} that are used in section 3.9. We plan
to carry out the details of this argument elsewhere.

\subsection*{Acknowledgments}
The authors would like to thank D. Dolgopyat, M. Brin and L.
Polterovich for useful comments regarding the proof of Lemma
\ref{dynlemma} and Yu. Safarov for suggesting the use of Riesz means
in the proof of Theorem \ref{Riesz}. We also thank Y. Colin de
Verdi\`{e}re, A.~Grigor'yan  and M. Shubin  for helpful discussions, and the
anonymous referee for useful remarks.

Part of this paper was written while the first author visited IHES
and Max Planck Institute for Mathematics in Bonn; their hospitality
is greatly appreciated.

\end{document}